\documentclass[12pt,twoside,a4paper]{article}
\usepackage[english]{babel}
\usepackage{amssymb}
\usepackage[T1]{fontenc}
\usepackage{amsmath,epsf}
\newtheorem{Th}{Theorem}
\newtheorem{lemma}{Lemma}

\newtheorem{Prop}{Proposition}

\newcommand{\I}{\text{I}}
\newcommand{\E}{\mathbb{E}}
\newcommand{\F}{\mathcal{F}}
\newcommand{\R}{\mathbb{R}}
\newcommand{\Z}{\mathbb{Z}}
\newcommand{\N}{\mathbb{N}}

\newcommand{\B}{\mathcal{B}}
\newcommand{\G}{\mathcal{G}}
\newcommand{\U}{\mathcal{U}}

\renewcommand{\P}{\mathbb{P}}

\newcounter{tictac}
\newenvironment{fleuveA}{
   \begin{list}{($\textbf{A\arabic{tictac}}$) }{\usecounter{tictac}
\leftmargin 1cm\labelwidth 2em}}{\end{list}}
\def\1{\,\rlap{\mbox{\small\rm 1}}\kern.15em 1}
\def\ind#1{\1_{#1}}
\def\build#1_#2^#3{\mathrel{\mathop{\kern 0pt#1}\limits_{#2}^{#3}}}
\def\tend#1#2{\build\hbox to 12mm{\rightarrowfill}_{#1\rightarrow #2}^{a.s.}}

\def\converge#1#2#3{\build\hbox to
15mm{\rightarrowfill}_{#1\rightarrow #2}^{\hbox{\scriptsize #3}}}
\begin{document}

\begin{center}
{\large
    {\sc
Kernel deconvolution estimation for random fields
    }
}
\bigskip

 Ahmed EL GHINI and Mohamed EL MACHKOURI

\end{center}

{\renewcommand\abstractname{Abstract}
\begin{abstract}
\baselineskip=18pt % espacement entre les lignes
In this work, we establish the asymptotic normality of the deconvolution kernel density estimator in the context of strongly mixing random fields. Only minimal 
conditions on the bandwidth parameter are required and a simple criterion on the strong mixing coefficients is provided. Our approach is based on the 
Lindeberg's method rather than on Bernstein's technique and coupling arguments widely used in previous works on nonparametric estimation for spatial processes. 
We deal also with nonmixing random fields which can be written as a (nonlinear) functional of i.i.d. random fields by considering the physical dependence measure coefficients 
introduced by Wu \cite{Wu2005}.
{\em Short title:} Deconvolution kernel density estimator for random fields.
\end{abstract}
%-------------------------------------------------- Introduction ----------------------------------------------------------
\thispagestyle{empty}
\baselineskip=18pt
\section{Introduction and main results}
Let $X=(X_i)_{i\in\Z^d}$ be a stationary real random field defined on the probability space $(\Omega,\F,\P)$. We observe the random field $X$ on a region $\Lambda_n$, $n\in\N^{\ast}$, 
but the observations are contamined with noise such as measurement errors. In fact, we observe only the random field $Y=(Y_i)_{i\in\Z^d}$ defined for any $i$ in $\Z^d$ by 
$Y_i=X_i+\theta_i$ where the error variables $(\theta_i)_{i\in\Z^d}$ are identically distributed and independent of $X$. 
We denote by $f_Y$, $f_X$ and $f_{\theta}$ the marginal density of $Y$, $X$ and $\theta$ respectively and we have $f_Y=f_X\star f_\theta$. 
We observe a sample of $Y$ and we want to estimate $f_X$ using the deconvolution kernel approach introduced by Stefanski and Carroll \cite{Stefanski--Carroll1990}. Previous 
key results on deconvolution kernel density estimators for time series are Fan \cite{Fan1991b}, \cite{Fan1991a} and Masry \cite{Masry1993b}, \cite{Masry1993a}. 
For strongly mixing random fields indexed by the lattice $\Z^d$, Li \cite{Li2008} obtained a central limit theorem for the deconvolution kernel density estimator using  
the so-called Bernstein's small and large blocks technique and coupling arguments initiated by Tran \cite{Tran1990}. Note that the extension 
of asymptotic result for time series to the spatial setting is not trivial because of difficulties coming from spatial ordering. The purpose of this work is to put on light 
a new approach for the asymptotic normality of kernel density estimators. In fact, we are going to apply the Lindeberg's method (see \cite{Lindeberg1922}) in 
order to improve the result by Li \cite{Li2008} in several directions. This new approach was recently applied successfully in 
El Machkouri and Stoica \cite{Elmachkouri--Stoica2010} and El Machkouri \cite{Elmachkouri2011} for the Nadaraya-Watson estimator and the Parzen-Rosenblatt estimator 
respectively in the setting of random fields.\\
For any finite subset $B$ of $\Z^d$, denote $\vert B\vert$ the number of elements in $B$ and $\partial B$ its boundary defined by 
$\partial B=\left\{i\in B\,;\,\exists j\notin B\,\,\vert i-j\vert=1\right\}$ where $\vert s\vert=\max_{1\leq k\leq d}\vert s_k\vert$ 
for any $s=(s_1,...,s_d)$ in $\Z^d$. In the sequel, we assume that we observe $(X_i)_{i\in\Z^d}$ on a sequence $(\Lambda_n)_{n\geq 1}$ 
of finite subsets of $\Z^d$ which satisfies 
\begin{equation}\label{condition_on_Lambda_n}
\lim_{n\to\infty}\vert\Lambda_n\vert=\infty\quad\textrm{and}\quad\lim_{n\to\infty}\frac{\vert\partial\Lambda_n\vert}{\vert\Lambda_n\vert}=0.
\end{equation}
Given two $\sigma$-algebras $\U$ and $\mathcal{V}$ of $\F$, we recall the $\alpha$-mixing coefficient introduced by Rosenblatt \cite{Ros} defined by 
$\alpha(\U,\mathcal{V})=\sup\{\vert\P(A\cap B)-\P(A)\P(B)\vert\, ,\,A\in\U,\,B\in\mathcal{V}\}$. 
For any $\tau$ in $\N^{\ast}\cup\{\infty\}$ and any positive integer $n$, we consider the mixing coefficient $\alpha_{1,\tau}(n)$ defined by
$$
\alpha_{1,\tau}(n)=\sup\,\{\alpha(\sigma(X_k),\F_{B}),\,k\in\Z^d,\, \vert B\vert\leq\tau,\,\rho(B,\{k\})\geq n\},
$$
where $\F_{B}=\sigma(X_i\,;\,i\in B)$ and the distance $\rho$ is defined for any subsets $B_1$ and $B_2$ of $\Z^d$ by 
$\rho(B_{1},B_{2})=\min\{\vert i-j\vert,\,i\in B_{1},\,j\in B_{2}\}$. 
We say that the random field $(X_i)_{i\in\Z^d}$ is strongly mixing if $\lim_{n\to\infty}\alpha_{1,\tau}(n)=0$ for some $\tau$ in $\N^{\ast}\cup\{\infty\}$.\\
Let $(b_n)_{n\geq 1}$ be a sequence of positive numbers going to zero as $n$ goes to infinity. The deconvolving kernel density estimator for $f_X$ is defined for 
any $x$ in $\R$ by 
\begin{equation}\label{definition_estimateur}
\hat{f}_n(x)=\frac{1}{\vert\Lambda_n\vert b_n}\sum_{i\in\Lambda_n}g_n\left(\frac{x-Y_i}{b_n}\right)
\end{equation}
where for any $z$ in $\R$,
$$
g_n(z)=\frac{1}{2\pi}\int_{\R} e^{-itz}\frac{\phi_K(t)}{\phi_{\theta}(t/b_n)}dt.
$$
The kernel density estimator $\hat{f}_n$ defined by $(\ref{definition_estimateur})$ can be written for any $x$ in $\R$ as
\begin{equation}\label{definition_estimateur_bis}
\hat{f}_n(x)=\frac{1}{2\pi}\int_{\R} e^{-itx}\hat{\phi}_n(t)\frac{\phi_K(tb_n)}{\phi_{\theta}(t)}dt
\end{equation}
where 
$$
\hat{\phi}_n(t)=\frac{1}{\vert\Lambda_n\vert}\sum_{i\in\Lambda_n}e^{itY_i}.
$$
We consider the following assumptions:
\begin{fleuveA}
\item The marginal probability distribution of each $Y_k$ is absolutely continuous with continuous positive density function $f_Y$.
\item There exists $\kappa>0$ such that $\sup_{(x,y)\in\R^2}f_{i\vert j}(y\vert x)\leq \kappa$ where $f_{i\vert j}$ is the conditional density 
function of $Y_i$ given $Y_j$ for any $i$ and $j$ in $\Z^d$.
\item There exists $\beta>0$ and $B>0$ such that $\vert t\vert^{\beta}\vert\phi_{\theta}(t)\vert\converge{t}{+\infty}{ }B$.
\item The characteristic function $\phi_K$ of the kernel $K$ vanishes outside $[-1,1]$.
\item The bandwidth $b_n$ converges to zero and $\vert\Lambda_n\vert  b_n$ goes to infinity.
\end{fleuveA}
The following result establishes the asymptotic bias of the estimator $\hat{f}_n$.
\begin{Prop}[Li \cite{Li2008}, 2008]\label{biais_asymptotique}
If $\phi_K$ is continuous then, for any real $x$, 
$$
\E(\hat{f}_n(x))\converge{n}{+\infty}{ }f(x).
$$ 
\end{Prop}
Now, we investigate the asymptotic variance of the estimator $\hat{f}_n$.
\begin{Prop}\label{variance_asymptotique}
Assume that $\sum_{m\geq 1}m^{2d-1}\,\alpha_{1,1}(m)<\infty$. For any $x$ in $\R$, we have
\begin{equation}\label{limite_variance_fn}
\lim_{n\to\infty}\vert\Lambda_n\vert  b_n^{2\beta+1}\mathbb{V}(\hat{f}_n(x))=\frac{f_Y(x)}{B^2}\int_{\R} \vert t\vert^{2\beta}\vert\phi_K(t)\vert^2dt:=\sigma^2(x).
\end{equation}
\end{Prop}
Our main result is the following.
\begin{Th}\label{convergence-loi}
Assume that Assumptions $\textbf{\emph{(A1)}},...,\textbf{\emph{(A5)}}$ hold and 
\begin{equation}\label{mixing-condition}
\sum_{m=1}^{+\infty}m^{2d-1}\,\alpha_{1,\infty}(m)<\infty.
\end{equation}
Then for any positive integer $k$ and any distinct points $x_1,...,x_k$ in $\R$,
\begin{equation}\label{limit}
(\vert\Lambda_n\vert  b_{n}^{2\beta+1})^{1/2}
\left(\begin{array}{c}
       \hat{f}_{n}(x_1)-\E \hat{f}_n(x_1)\\
       \vdots\\
       \hat{f}_{n}(x_k)-\E\hat{f}_n(x_k)
       \end{array} \right)
\converge{n}{+\infty}{\textrm{$\mathcal{L}$}}
\mathcal{N}\left(0,V\right)
\end{equation}
where $V$ is a diagonal matrix with diagonal elements $v_{ii}=\frac{f_Y(x_i)}{B^2}\int_{\R} \vert t\vert^{2\beta}\vert\phi_K(t)\vert^2dt$.
\end{Th}
\textbf{Remark 1}. Theorem \ref{convergence-loi} improves Theorem 4.1 and Theorem 4.2 in \cite{Li2008} in three directions: the regions $\Lambda_n$ where 
the random field is observed are not reduced to rectangular ones, the assumption $\textbf{(A5)}$ on the bandwidth $b_n$ is minimal and the mixing condition 
$(\ref{mixing-condition})$ does not dependent on the bandwith parameter $b_n$.\\
\\
Since the mixing property is often unverifiable and might be too restrictive, it is important to provide limit theorems for 
nonmixing and possibly nonlinear spatial processes. So, in the sequel, we consider that $(X_i)_{i\in\Z^d}$ is a field of 
identically distributed real random variables with a marginal density $f$ such that 
\begin{equation}\label{definition_champ}
X_i = F\left(\varepsilon_{i-s};\,s\in\Z^d\right),\quad i\in\Z^d,
\end{equation}
where $(\varepsilon_j)_{j\in\Z^d}$ are i.i.d. random variables and $F$ is a measurable function defined on $\R^{\Z^d}$. 
In the one-dimensional case ($d=1$), the class (\ref{definition_champ}) includes linear as well as many widely used nonlinear time series models as special cases. More importantly, 
it provides a very general framework for asymptotic theory for statistics of stationary time series (see \cite{Wu2005} and the review paper \cite{Wu2011}). 
Let $(\varepsilon_j^{'})_{j\in\Z^d}$ be an i.i.d. copy of $(\varepsilon_j)_{j\in\Z^d}$ and consider for any positive integer $n$ the coupled version 
$X_i^{\ast}$ of $X_i$ defined by $X_i^{\ast}=F\left(\varepsilon^{\ast}_{i-s}\,;\,s\in\Z^d\right)$ where 
$\varepsilon_j^{\ast}=\varepsilon_j\ind{\{j\neq 0\}}+\varepsilon_0^{'}\ind{\{j=0\}}$ for any $j$ in $\Z^d$. In other words, we obtain $X_i^{\ast}$ from $X_i$ 
by just replacing $\varepsilon_0$ by its copy $\varepsilon_0^{'}$. Following Wu \cite{Wu2005}, 
we introduce appropriate dependence measures: let $i$ in $\Z^d$ and $p>0$ be fixed. If $X_i$ belongs to $\mathbb{L}_{p}$ (that is, $\E\vert X_i\vert^p$ is finite), 
we define the physical dependence measure $\delta_{i,p}=\|X_i-X_i^{\ast}\|_{p}$ where $\|\,.\,\|_p$ is the usual $\mathbb{L}^p$-norm and we say that the random field 
$(X_i)_{i\in\Z^d}$ is $p$-stable if $\sum_{i\in\Z^d}\delta_{i,p}<\infty$. For $d\geq 2$, the reader should keep in mind the following two examples already given in \cite{Elmachkouri-Volny-Wu2011} :\\
\underline{{\em Linear random fields}}: Let $(\varepsilon_i)_{i\in\Z^d}$ be
i.i.d random variables with $\varepsilon_i$ in $\mathbb{L}^p$, $p
\geq 2$. The linear random field $X$ defined for any $i$ in $\Z^d$
by
$$
X_i=\sum_{s\in\Z^d}a_s\varepsilon_{i-s}
$$
with $(a_s)_{s\in\Z^d}$ in $\R^{\Z^d}$ such that $\sum_{i\in\Z^d}a_i^2<\infty$ is of the form $(\ref{definition_champ})$ with a linear functional
$g$. For any $i$ in $\Z^d$, $\delta_{i,p} = \|a_i\| \|
\varepsilon_0 - \varepsilon^{'}_0\|_p$. So, $X$ is $p$-stable if $\sum_{i\in\Z^d}\vert a_i\vert<\infty$. 
Clearly, if $\textrm{H}$ is a Lipschitz continuous function, under the
above condition, the subordinated process $Y_i = \textrm{H}(X_i)$ is also $p$-stable since $\delta_{i,p} = O(|a_i|)$. \\
\underline{{\em Volterra field}} : Another class of nonlinear random field is
the Volterra process which plays an important role in the
nonlinear system theory (Casti \cite{Casti1985}, Rugh \cite{Rugh1981}): consider the
second order Volterra process
\begin{eqnarray*}
X_i = \sum_{s_1, s_2\in\Z^d} a_{s_1, s_2} \varepsilon_{i-s_1} \varepsilon_{i-s_2},
\end{eqnarray*}
where $a_{s_1, s_2}$ are real coefficients with $a_{s_1, s_2} = 0$ if $s_1 = s_2$ and $(\varepsilon_i)_{i\in\Z^d}$ are i.i.d. random variables with $\varepsilon_i$ 
in $\mathbb{L}^p$, $p \geq 2$.
Let
\begin{eqnarray*}
A_i = \sum_{s_1, s_2\in\Z^d} (a_{s_1, i}^2 + a_{i, s_2}^2)\quad\textrm{and}\quad 
B_i = \sum_{s_1, s_2\in\Z^d} (|a_{s_1, i}|^p + |a_{i, s_2}|^p).
\end{eqnarray*}
By the Rosenthal inequality, there exists a constant $C_p > 0$ such that
\begin{eqnarray*}
\delta_{i,p} = \| X_i - X_i^*\|_p \leq C_p A_i^{1/2} \| \varepsilon_0\|_2 \| \varepsilon_0\|_p+ C_p B_i^{1/p} \| \varepsilon_0\|_p^2.
\end{eqnarray*}
\begin{Th}\label{convergence-loi-nonmixing-case}
Let $(X_i)_{i\in\Z^d}$ be defined by the relation ($\ref{definition_champ}$) and assume that $\textbf{\emph{(A1)}},...,\textbf{\emph{(A5)}}$ hold. 
If ($\ref{mixing-condition}$) is replaced by the condition  
\begin{equation}\label{nonmixing-condition}
\sum_{i\in\Z^d}\vert i\vert^{\frac{5d}{2}}\delta_i<\infty
\end{equation}
then the conclusion of Theorem \ref{convergence-loi} still holds.  
\end{Th}
\section{Proofs}
Throughout this section, the symbol $\kappa$ will denote a generic positive constant which the value is not important and 
for any $i=(i_1,...,i_d)\in\Z^d$, we denote $\vert i\vert=\max_{1\leq k\leq d}\vert i_k\vert$. Recall also that for any finite subset $B$ of $\Z^d$, 
we denote $\vert B\vert$ the number of elements in $B$. Let $\tau\in\N^{\ast}\cup\{\infty\}$ be fixed and consider the sequence $(m_{n,\tau})_{n\geq 1}$ defined by 
\begin{equation}\label{def_mn_tau}
m_{n,\tau}=\max\left\{v_n,\left[\left(\frac{1}{b_n}\sum_{\vert i\vert>v_n}\vert i\vert^d\,\alpha_{1,\tau}(\vert i\vert)\right)^{\frac{1}{d}}\right]+1\right\}
\end{equation}
where $v_n=\big[b_n^{\frac{-1}{2d}}\big]$ and $[\,.\,]$ denotes the integer part function. The following technical lemma is a spatial version of a 
result by Bosq, Merlev\`ede and Peligrad (\cite{Bosq-Merlevede-Peligrad1999}, pages 88-89).
\begin{lemma}\label{mn-tau}
If $\tau\in\N^{\ast}\cup\{\infty\}$ and $\sum_{m\geq 1}m^{2d-1}\,\alpha_{1,\tau}(m)<\infty$ then
$$
m_{n,\tau}\to\infty,\quad m_{n,\tau}^d\,b_n\to0\quad\textrm{and}\quad\frac{1}{m_{n,\tau}^d\,b_n}\sum_{\vert i\vert>m_{n,\tau}}\vert i\vert^d\,\alpha_{1,\tau}(\vert i\vert)\to 0.
$$
\end{lemma}
{\em Proof of Proposition \ref{variance_asymptotique}}. Let $z$ be fixed in $\R$. For any $i$ in $\Z^d$, we denote 
$$
Z_i(z)=\frac{1}{b_n}\left(g_n\left(\frac{z-Y_i}{b_n}\right)-\E g_n\left(\frac{z-Y_i}{b_n}\right)\right).
$$
The proof of the following lemma is done in the appendix.
\begin{lemma}\label{lemme_Z}
For any $z$ in $\R$,
\begin{equation}\label{convergence_esperance_Zs_carre}
b_n^{2\beta+1}\E(Z_0^2(z))\converge{n}{\infty}{ }\sigma^2(z):=\frac{f_Y(z)}{B^2}\int_{\R} \vert u\vert^{2\beta}\vert\phi_K(u)\vert^2du
\end{equation}
and $\sup_{i\in\Z^d\backslash\{0\}}\E\left\vert Z_0(s)Z_i(t)\right\vert=O(b_n^{-2\beta})$ for any $s$ and $t$ in $\R$.
\end{lemma}
Let $x$ in $\R$ be fixed. We have 
\begin{equation}\label{decomposition_variance_fnx}
\vert\Lambda_n\vert  b_n^{2\beta+1}\mathbb{V}(\hat{f}_n(x))=b_n^{2\beta+1}\E\left(Z_0^2(x)\right)
+\frac{b_n^{2\beta+1}}{\vert\Lambda_n\vert  }\sum_{\stackrel{i,j\in\Lambda_n}{i\neq j}}\textrm{Cov}(Z_i(x),Z_j(x)).
\end{equation}
Since $(Z_i)_{i\in\Z^d}$ is stationary, we have 
\begin{align*}
\frac{b_n^{2\beta+1}}{\vert\Lambda_n\vert  }\left\vert\sum_{\stackrel{i,j\in\Lambda_n}{i\neq j}}\textrm{Cov}(Z_i(x),Z_j(x))\right\vert
&\leq b_n^{2\beta+1}\sum_{i\in\Z^d\backslash\{0\}}\vert\E\left(Z_0(x)Z_i(x)\right)\vert\\
&\leq b_n^{2\beta+1}\left(m_{n,1}^d\sup_{i\in\Z^d\backslash\{0\}}\E\left\vert Z_0(x)Z_i(x)\right\vert
+\sum_{\vert i\vert>m_{n,1}}\vert\E\left(Z_0(x)Z_i(x)\right)\vert\right).
\end{align*}
By Rio's covariance inequality (cf. \cite{Rio1993}, Theorem 1.1), we know that $\vert\E\left(Z_0(x)Z_i(x)\right)\vert\leq \kappa \|Z_0(x)\|_{\infty}^2\alpha_{1,1}(\vert i\vert)$. 
Since $\|Z_0(x)\|_{\infty}\leq \kappa b_n^{-\beta-1}$ and $\tau\in\N^{\ast}\cup\{\infty\}$, we obtain 
$$
\frac{b_n^{2\beta+1}}{\vert\Lambda_n\vert  }\left\vert\sum_{\stackrel{i,j\in\Lambda_n}{i\neq j}}\textrm{Cov}(Z_i(x),Z_j(x))\right\vert
\leq\kappa\left(m_{n,1}^db_n^{2\beta+1}\sup_{i\in\Z^d\backslash\{0\}}\E\left\vert Z_0(x)Z_i(x)\right\vert
+\frac{1}{m_{n,1}^db_n}\sum_{\vert i\vert>m_{n,1}}\vert i\vert^d\alpha_{1,1}(\vert i\vert)\right).
$$
Applying Lemma \ref{mn-tau} and the second part of Lemma \ref{lemme_Z}, we derive 
\begin{equation}\label{limite_covariance_Zix_Zjx}
\lim_{n\to\infty}\frac{b_n^{2\beta+1}}{\vert\Lambda_n\vert  }\left\vert\sum_{\stackrel{i,j\in\Lambda_n}{i\neq j}}\textrm{Cov}(Z_i(x),Z_j(x))\right\vert=0.
\end{equation}
Combining (\ref{convergence_esperance_Zs_carre}), (\ref{decomposition_variance_fnx}) and (\ref{limite_covariance_Zix_Zjx}), we obtain 
($\ref{limite_variance_fn}$). The proof of Proposition \ref{variance_asymptotique} is complete.\\
\\
\textbf{Proof of Theorem $\textbf{\ref{convergence-loi}}$}. Without loss of generality, we consider only the case $k=2$ and we 
refer to $x_1$ and $x_2$ as $x$ and $y$ ($x\neq y$). Let $\lambda_1$ and $\lambda_2$ be two constants such that $\lambda_1^2+\lambda_2^2=1$ and 
denote 
$$
S_n=\lambda_1(\vert\Lambda_n\vert  b_n^{2\beta+1})^{1/2}(\hat{f}_n(x)-\E\hat{f}_n(x))+\lambda_2(\vert\Lambda_n\vert  b_n^{2\beta+1})^{1/2}(\hat{f}_n(y)-\E\hat{f}_n(y))
=\sum_{i\in\Lambda_n}\frac{b_n^{\beta+\frac{1}{2}}\Delta_i}{\vert\Lambda_n\vert^{1/2}}
$$ 
where $\Delta_i=\lambda_1Z_{i}(x)+\lambda_2Z_{i}(y)$ and for any $z$ in $\R$,
$$
Z_{i}(z)=\frac{1}{b_n}\left(g_n\left(\frac{z-Y_{i}}{b_n}\right)-\E g_n\left(\frac{z-Y_{i}}{b_n}\right)\right).
$$
We consider the notation
\begin{equation}\label{eta_sigma}
\eta=\frac{\lambda_1^2f_Y(x)+\lambda_2^2f_Y(y)}{B^2}\int_{\R} \vert t\vert^{2\beta}\vert\phi_K(t)\vert^2dt.
\end{equation}
The proof of the following technical result is postponed to the annex.
\begin{lemma}\label{lemme-technique} $b_n^{2\beta+1}\E(\Delta_0^2)$ converges to $\eta$ as $n$ goes to infinity and 
$\sup_{i\in\Z^d\backslash\{0\}}\E\vert\Delta_0\Delta_i\vert=O(b_n^{-2\beta})$.
\end{lemma}
In order to prove the convergence in distribution of $S_n$ to $\sqrt{\eta}\tau_0$ where $\tau_0\sim\mathcal{N}(0,1)$, 
we follow the Lindeberg's method used in the proof of the central limit theorem for stationary random 
fields by Dedecker \cite{Dedecker1998}. Let $\varphi$ be a one to one map from $[1,\kappa]\cap\N^{\ast}$ to a finite
subset of $\Z^d$ and $(\xi_i)_{i\in\Z^d}$ a real random field. For all integers $k$ in $[1,\kappa]$, we denote
$$
S_{\varphi(k)}(\xi)=\sum_{i=1}^k \xi_{\varphi(i)}\quad\textrm{and}\quad
S_{\varphi(k)}^{c}(\xi)=\sum_{i=k}^\kappa \xi_{\varphi(i)}
$$
with the convention $S_{\varphi(0)}(\xi)=S_{\varphi(\kappa+1)}^{c}(\xi)=0$. To describe the set $\Lambda_{n}$, we define the one to
one map $\varphi$ from $[1,\vert\Lambda_n\vert  ]\cap\N^{\ast}$ to $\Lambda_{n}$ by: $\varphi$ is the unique function such that 
$\varphi(k)<_{\text{\text{lex}}}\varphi(l)$ for $1\leq k<l\leq\vert\Lambda_n\vert  $. From now on, we consider a field $(\tau_{i})_{i\in\Z^d}$ of i.i.d. 
random variables independent of $(\Delta_{i})_{i\in\Z^d}$ such that $\tau_{0}$ has the standard normal law $\mathcal{N}(0,1)$. 
We introduce the fields $\Gamma$ and $\gamma$ defined for any $i$ in $\Z^d$ by 
$$
\Gamma_{i}=\frac{b_n^{\beta+\frac{1}{2}}\Delta_i}{\vert\Lambda_n\vert^{1/2}}\quad\textrm{and}\quad\gamma_{i}=\frac{\tau_i\sqrt{\eta}}{\vert\Lambda_n\vert^{1/2}}
$$ 
where $\eta$ is defined by ($\ref{eta_sigma}$). Let $h$ be any function from $\R$ to $\R$. For $0\leq k\leq l\leq \vert\Lambda_n\vert+1$, 
we introduce $h_{k,l}(\Gamma)=h(S_{\varphi(k)}(\Gamma)+S_{\varphi(l)}^{c}(\gamma))$. 
With the above convention we have that $h_{k,\vert\Lambda_n\vert  +1}(\Gamma)=h(S_{\varphi(k)}(\Gamma))$ and also $h_{0,l}(\Gamma)=h(S_{\varphi(l)}^{c}(\gamma))$. 
In the sequel, we will often write $h_{k,l}$ instead of $h_{k,l}(\Gamma)$. We denote by $B_{1}^4(\R)$ the unit ball of $C_{b}^4(\R)$: $h$ belongs to
$B_{1}^4(\R)$ if and only if it belongs to $C^4(\R)$ and satisfies $\max_{0\leq i\leq 4}\|h^{(i)}\|_{\infty}\leq 1$. It suffices to prove that for all $h$ in $B_{1}^4(\R)$,
$$
\E\left(h\left(S_{\varphi(\vert\Lambda_n\vert  )}(\Gamma)\right)\right)\converge{n}{+\infty}{}\E \left(h\left(\tau_0\sqrt{\eta}\right)\right).
$$
We use Lindeberg's decomposition:
$$
\E \left(h\left(S_{\varphi(\vert\Lambda_n\vert  )}(\Gamma)\right)-h\left(\tau_{0}\sqrt{\eta}\right)\right)=\sum_{k=1}^{\vert\Lambda_n\vert  }\E \left(h_{k,k+1}-h_{k-1,k}\right).
$$
Now,
$$
h_{k,k+1}-h_{k-1,k}=h_{k,k+1}-h_{k-1,k+1}+h_{k-1,k+1}-h_{k-1,k}.
$$
Applying Taylor's formula we get that:
$$
h_{k,k+1}-h_{k-1,k+1}=\Gamma_{\varphi(k)}h_{k-1,k+1}^{'}+\frac{1}{2}\Gamma_{\varphi(k)}^{2}h_{k-1,k+1}^{''}+R_{k}
$$
and
$$
h_{k-1,k+1}-h_{k-1,k}=-\gamma_{\varphi(k)}h_{k-1,k+1}^{'}-\frac{1}{2}\gamma_{\varphi(k)}^{2}h_{k-1,k+1}^{''}+r_{k}
$$
where $\vert R_{k}\vert\leq \Gamma_{\varphi(k)}^2(1\wedge\vert \Gamma_{\varphi(k)}\vert)$ and 
$\vert r_{k}\vert\leq\gamma_{\varphi(k)}^2(1\wedge\vert\gamma_{\varphi(k)}\vert)$.
Since $(\Gamma,\tau_{i})_{i\neq \varphi(k)}$ is independent of $\tau_{\varphi(k)}$, it follows that
$$
\E \left(\gamma_{\varphi(k)}h_{k-1,k+1}^{'}\right)=0\quad\textrm{and}\quad
\E \left(\gamma_{\varphi(k)}^2h_{k-1,k+1}^{''}\right)=\E \left(\frac{\eta}{\vert\Lambda_n\vert  }h_{k-1,k+1}^{''}\right)
$$
Hence, we obtain
\begin{align*}
\E \left(h(S_{\varphi(\vert\Lambda_n\vert  )}(\Gamma))-h\left(\tau_0\sqrt{\eta}\right)\right)&=\sum_{k=1}^{\vert\Lambda_n\vert  }\E (\Gamma_{\varphi(k)}h_{k-1,k+1}^{'})\\
&\quad+\sum_{k=1}^{\vert\Lambda_n\vert  }\E \left(\left(\Gamma_{\varphi(k)}^2-\frac{\eta}{\vert\Lambda_n\vert  }\right)\frac{h_{k-1,k+1}^{''}}{2}\right)\\
&\quad+\sum_{k=1}^{\vert\Lambda_n\vert  }\E \left(R_{k}+r_{k}\right).
\end{align*}
Let $1\leq k\leq \vert\Lambda_n\vert  $ be fixed. Noting that $\vert\Delta_0\vert$ is bounded by $\kappa b_n^{-\beta-1}$ and applying the first part of 
Lemma $\ref{lemme-technique}$, we derive
$$
\E\vert R_k\vert\leq\frac{b_n^{3\beta+\frac{3}{2}}\E\vert\Delta_0\vert^3}{\vert\Lambda_n\vert^{3/2}}=O\left(\frac{1}{(\vert\Lambda_n\vert^{3}\,b_n)^{1/2}}\right)
$$
and
$$
\E\vert r_k\vert\leq\frac{\E\vert\gamma_0\vert^3}{\vert\Lambda_n\vert^{3/2}}\leq\frac{\eta^{3/2}\E\vert\tau_0\vert^3}{\vert\Lambda_n\vert^{3/2}}
=O\left(\frac{1}{\vert\Lambda_n\vert^{3/2}}\right).
$$
Consequently, we obtain
$$
\sum_{k=1}^{\vert\Lambda_n\vert}\E\left(\vert R_{k}\vert+\vert r_{k}\vert\right)
=O\left(\frac{1}{(\vert\Lambda_n\vert  b_n)^{1/2}}+\frac{1}{\vert\Lambda_n\vert^{1/2}}\right)=o(1).
$$
Now, it is sufficient to show
\begin{equation}\label{equation1_mixing_case}
\lim_{n\to+\infty}\sum_{k=1}^{\vert\Lambda_n\vert  }\left(\E (\Gamma_{\varphi(k)}h_{k-1,k+1}^{'})+\E \left(\left(\Gamma_{\varphi(k)}^2-\frac{\eta}{\vert\Lambda_n\vert  }\right)\frac{h_{k-1,k+1}^{''}}{2}\right)\right)=0.
\end{equation}
First, we focus on
$\sum_{k=1}^{\vert\Lambda_n\vert  }\E \left(\Gamma_{\varphi(k)}h_{k-1,k+1}^{'}\right)$. On the lattice $\Z^{d}$ we define the
lexicographic order as follows: if $i=(i_{1},...,i_{d})$ and $j=(j_{1},...,j_{d})$ are distinct elements of $\Z^{d}$, the
notation $i<_{\text{\text{lex}}}j$ means that either $i_{1}<j_{1}$ or for some $p$ in $\{2,3,...,d\}$, $i_{p}<j_{p}$ and $i_{q}=j_{q}$ for $1\leq
q<p$. Let the sets $\{V_{i}^{M}\,;\,i\in\Z^{d}\,,\,M\in\N^{\ast}\}$ be defined as follows:
$$
V_{i}^{1}=\{j\in\Z^{d}\,;\,j<_{\text{lex}}i\}\,\,\,\textrm{and for $M\geq 2$},\,\,V_{i}^{M}=V_{i}^{1}\cap\{j\in\Z^{d}\,;\,\vert i-j\vert\geq M\}.
$$
For any subset $L$ of $\Z^{d}$ define $\F_{L}=\sigma(\Delta_{i}\,;\,i\in L)$ and set
$$
\E_{M}(\Delta_{i})=\E(\Delta_{i}\vert\F_{V_{i}^{M}}),\quad M\in \N^{\ast}.
$$
For all $M$ in $\N^{\ast}$ and all integer $k$ in $[1,\vert\Lambda_n\vert  ]$, we define
$$
\textrm{E}_{k}^M=\varphi([1,k]\cap\N^{\ast})\cap V_{\varphi(k)}^M\quad\textrm{and}\quad S_{\varphi(k)}^M(\Gamma)=\sum_{i\in \textrm{E}_{k}^M}\Gamma_{i}.
$$
For any function $\Psi$ from $\R$ to $\R$, we define $\Psi_{k-1,l}^M=\Psi(S_{\varphi(k)}^M(\Gamma)+S_{\varphi(l)}^c(\gamma))$. We are going to 
apply this notation to the successive derivatives of the function $h$. Our aim is to show that
$$
\lim_{n\to+\infty}\sum_{k=1}^{\vert\Lambda_n\vert  }\E\left(\Gamma_{\varphi(k)}h_{k-1,k+1}^{'}-\Gamma_{\varphi(k)}\left(S_{\varphi(k-1)}(\Gamma)-S_{\varphi(k)}^{m_{n,\infty}}(\Gamma)\right)h_{k-1,k+1}^{''}\right)=0.
$$
where $(m_{n,\infty})_{n\geq 1}$ is the sequence defined by $(\ref{def_mn_tau})$. First, we use the decomposition
$$
\Gamma_{\varphi(k)}h_{k-1,k+1}^{'}=\Gamma_{\varphi(k)}h_{k-1,k+1}^{'m_{n,\infty}}+\Gamma_{\varphi(k)}\left(h_{k-1,k+1}^{'}-h_{k-1,k+1}^{'m_{n,\infty}}\right).
$$
We consider a one to one map $\psi$ from $[1,\vert \textrm{E}_{k}^{m_{n,\infty}}\vert]\cap\N^{\ast}$ to $\textrm{E}_{k}^{m_{n,\infty}}$ and such that $\vert
\psi(i)-\varphi(k)\vert\leq\vert \psi(i-1)-\varphi(k)\vert$. This choice of $\psi$ ensures that $S_{\psi(i)}(\Gamma)$ and $S_{\psi(i-1)}(\Gamma)$ are
$\F_{V_{\varphi(k)}^{\vert \psi(i)-\varphi(k)\vert}}$-measurable. The fact that $\gamma$ is independent of $\Gamma$ imply that
$$
\E \left(\Gamma_{\varphi(k)}h^{'}\left(S_{\varphi(k+1)}^c(\gamma)\right)\right)=0.
$$
Therefore
\begin{equation}\label{equation_theta}
\left\vert \E \left(\Gamma_{\varphi(k)}h_{k-1,k+1}^{'{m_{n,\infty}}}\right)\right\vert=\left\vert\sum_{i=1}^{\vert
\textrm{E}_{k}^{m_{n,\infty}}\vert}\E \left(\Gamma_{\varphi(k)}\left(\Theta_i-\Theta_{i-1}\right)\right)\right\vert
\end{equation}
where $\Theta_i=h^{'}\left(S_{\psi(i)}(\Gamma)+S_{\varphi(k+1)}^c(\gamma)\right)$. 
Since $S_{\psi(i)}(\Gamma)$ and $S_{\psi(i-1)}(\Gamma)$ are $\F_{V_{\varphi(k)}^{\vert\psi(i)-\varphi(k)\vert}}$-measurable, 
we can take the conditional expectation of $\Gamma_{\varphi(k)}$ with respect to $\F_{V_{\varphi(k)}^{\vert\psi(i)-\varphi(k)\vert}}$ 
in the right hand side of ($\ref{equation_theta}$). On the other hand the function $h^{'}$ is $1$-Lipschitz, hence 
$\left\vert\Theta_i-\Theta_{i-1}\right\vert\leq\vert \Gamma_{\psi(i)}\vert$. Consequently,
$$
\left\vert\E \left(\Gamma_{\varphi(k)}\left(\Theta_i-\Theta_{i-1}\right)\right)\right\vert\leq\E \vert \Gamma_{\psi(i)}\E _{\vert \psi(i)-\varphi(k)\vert}\left(\Gamma_{\varphi(k)}\right)\vert
$$
and
$$
\left\vert
\E \left(\Gamma_{\varphi(k)}h_{k-1,k+1}^{'m_{n,\infty}}\right)\right\vert\leq\sum_{i=1}^{\vert
\textrm{E}_{k}^{m_{n,\infty}}\vert}\E \vert \Gamma_{\psi(i)}\E _{\vert \psi(i)-\varphi(k)\vert}(\Gamma_{\varphi(k)})\vert.
$$
Hence,
\begin{align*}
\left\vert\sum_{k=1}^{\vert\Lambda_n\vert  }\E \left(\Gamma_{\varphi(k)}h_{k-1,k+1}^{'m_{n,\infty}}\right)\right\vert
&\leq\frac{b_n^{2\beta+1}}{\vert\Lambda_n\vert  }\sum_{k=1}^{\vert\Lambda_n\vert  }\sum_{i=1}^{\vert \textrm{E}_{k}^{m_{n,\infty}}\vert} \E \vert\Delta_{\psi(i)}\E _{\vert \psi(i)-\varphi(k)\vert}(\Delta_{\varphi(k)})\vert\\
&\leq b_n^{2\beta+1}\sum_{\vert j\vert \geq m_{n,\infty}}\|\Delta_j\E _{\vert j\vert}(\Delta_0)\|_1.
\end{align*}
For any $j$ in $\Z^d$, we have 
$$
\|\Delta_j\E _{\vert j\vert}(\Delta_0)\|_1
=\textrm{Cov}\left(\vert\Delta_j\vert\left(\ind{\E_{\vert j\vert}(\Delta_0)\geq 0}-\ind{\E_{\vert j\vert}(\Delta_0)<0}\right),\Delta_0\right).
$$
So, applying Rio's covariance inequality (cf. \cite{Rio1993}, Theorem 1.1), we obtain
$$
\|\Delta_j\E _{\vert j\vert}(\Delta_0)\|_1\leq 4\int_{0}^{\alpha_{1,\infty}(\vert j\vert)}Q_{\Delta_0}^2(u)du
$$
where $Q_{\Delta_0}$ is defined by $Q_{\Delta_0}(u)=\inf\{t\geq 0\,;\,\P(\vert \Delta_0\vert>t)\leq u\}$ for any $u$ in $[0,1]$. Since $\vert\Delta_0\vert$ is bounded by 
$\kappa b_n^{-\beta-1}$, we have
$$
Q_{\Delta_0}(u)\leq\kappa b_n^{-\beta-1}\qquad\textrm{and}\qquad\|\Delta_j\E _{\vert j\vert}(\Delta_0)\|_1\leq \kappa b_n^{-2\beta-2}\alpha_{1,\infty}(\vert j\vert).
$$
Finally, we derive
$$
\left\vert\sum_{k=1}^{\vert\Lambda_n\vert  }\E \left(\Gamma_{\varphi(k)}h_{k-1,k+1}^{'m_{n,\infty}}\right)\right\vert
\leq \frac{\kappa}{b_n}\sum_{\vert j\vert \geq m_{n,\infty}}\alpha_{1,\infty}(\vert j\vert)
\leq\frac{\kappa}{m_{n,\infty}^db_n}\sum_{\vert j\vert \geq m_{n,\infty}}\vert j\vert^d\,\alpha_{1,\infty}(\vert j\vert)
$$
and by Lemma \ref{mn-tau}, we obtain $\lim_{n\to+\infty}\left\vert\sum_{k=1}^{\vert\Lambda_n\vert  }\E \left(\Gamma_{\varphi(k)}h_{k-1,k+1}^{'m_{n,\infty}}\right)\right\vert=0$. 
Applying again Taylor's formula, it remains to consider
$$
\Gamma_{\varphi(k)}(h_{k-1,k+1}^{'}-h_{k-1,k+1}^{'m_{n,\infty}})=\Gamma_{\varphi(k)}(S_{\varphi(k-1)}(\Gamma)-S_{\varphi(k)}^{m_{n,\infty}}(\Gamma))h_{k-1,k+1}^{''}+R_{k}^{'},
$$
where $\vert R_{k}^{'}\vert\leq 2\vert \Gamma_{\varphi(k)}(S_{\varphi(k-1)}(\Gamma)-S_{\varphi(k)}^{m_{n,\infty}}(\Gamma))
(1\wedge\vert S_{\varphi(k-1)}(\Gamma)-S_{\varphi(k)}^{m_{n,\infty}}(\Gamma)\vert)\vert$. Denoting $W_n=\{-m_{n,\infty}+1,...,m_{n,\infty}-1\}^d$ 
and $W_n^{\ast}=W_n\backslash\{0\}$, it follows that
\begin{align*}
\sum_{k=1}^{\vert\Lambda_n\vert  }\E \vert R_{k}^{'}\vert 
&\leq 2b_n^{2\beta+1}\E \left(\vert\Delta_{0}\vert\left(\sum_{i\in W_n}\vert\Delta_{i}\vert\right)
\left(1\wedge\frac{b_n^{\beta+\frac{1}{2}}}{\vert\Lambda_n\vert^{1/2}}\sum_{i\in W_n}\vert\Delta_{i}\vert\right)\right)\\
&=2b_n^{2\beta+1}\E\left(\left(\Delta_0^2+\sum_{i\in W_n^{\ast}}\vert\Delta_0\Delta_i\vert\right)\left(1\wedge\frac{b_n^{\beta+\frac{1}{2}}}{\vert\Lambda_n\vert^{1/2}}\sum_{i\in W_n}\vert\Delta_i\vert\right)\right)\\
&\leq\frac{2b_n^{3\beta+\frac{3}{2}}}{\vert\Lambda_n\vert^{1/2}}\sum_{i\in W_n}\E(\Delta_0^2\vert\Delta_i\vert)+2b_n^{2\beta+1}\sum_{i\in W_n^{\ast}}\E\vert\Delta_0\Delta_i\vert.
\end{align*}
Since $\vert\overline{\Delta}_0\vert$ is bounded by $\kappa b_n^{-\beta-1}$, we derive
\begin{align*}
\sum_{k=1}^{\vert\Lambda_n\vert  }\E \vert R_{k}^{'}\vert
&\leq\frac{\kappa b_n^{2\beta+\frac{1}{2}}}{\vert\Lambda_n\vert^{1/2}}\sum_{i\in W_n}\E(\vert\Delta_0\Delta_i\vert)+2b_n^{2\beta+1}\sum_{i\in W_n^{\ast}}\E\vert\Delta_0\Delta_i\vert\\
&=\frac{\kappa b_n^{2\beta+1}\E (\Delta_0^2)}{(\vert\Lambda_n\vert b_n)^{1/2}}+\left(\frac{\kappa b_n^{2\beta+1}}{(\vert\Lambda_n\vert b_n)^{1/2}}+2b_n^{2\beta+1}\right)\sum_{i\in W_n^{\ast}}\E(\vert\Delta_0\Delta_i\vert)\\
&\leq \kappa\left(\frac{1}{(\vert\Lambda_n\vert b_n)^{1/2}}+\left(\frac{b_n^{2\beta+1}}{(\vert\Lambda_n\vert b_n)^{1/2}}+b_n^{2\beta+1}\right)m_{n,\infty}^db_n^{-2\beta}\right)\qquad\textrm{(by Lemma $\ref{lemme-technique}$)}\\
&=\kappa\left(\frac{1}{(\vert\Lambda_n\vert b_n)^{1/2}}+\left(\frac{1}{(\vert\Lambda_n\vert b_n)^{1/2}}+1\right)m_{n,\infty}^db_n\right)\\
&=o(1)\qquad\textrm{(by Lemma $\ref{mn-tau}$ and Assumption $\textbf{(A5)}$}).
\end{align*}
So, we have shown that 
$$
\lim_{n\to +\infty}\sum_{k=1}^{\vert\Lambda_n\vert  }\E \left(\Gamma_{\varphi(k)}h^{'}_{k-1,k+1}-\Gamma_{\varphi(k)}(S_{\varphi(k-1)}(\Gamma)-S_{\varphi(k)}^{m_{n,\infty}}(\Gamma))h^{''}_{k-1,k+1}\right)=0.
$$
In order to obtain ($\ref{equation1_mixing_case}$) it remains to control
$$
F_{0}=\E \left(\sum_{k=1}^{\vert\Lambda_n\vert}h_{k-1,k+1}^{''}\left(\frac{\Gamma_{\varphi(k)}^2}{2}
+\Gamma_{\varphi(k)}\left(S_{\varphi(k-1)}(\Gamma)-S_{\varphi(k)}^{m_{n,\infty}}(\Gamma)\right)-\frac{\eta}{2\vert\Lambda_n\vert  }\right)\right).
$$
Let $\mu$ be the law of the stationary 
real random field $(\Delta_{i})_{i\in\Z^{d}}$ and consider the projection $\pi_0$ from $\R^{\Z^{d}}$ to $\R$ defined by $\pi_0(\omega)=\omega_{0}$ and
the family of translation operators $(T^{i})_{i\in\Z^{d}}$ from $\R^{\Z^{d}}$ to $\R^{\Z^{d}}$ defined by
$(T^{i}(\omega))_{k}=\omega_{k+i}$ for any $i\in\Z^{d}$ and any $\omega$ in $\R^{\Z^{d}}$. Denote by $\B$ the Borel
$\sigma$-algebra of $\R$. The random field $(\pi_0\circ T^{i})_{i\in\Z^{d}}$ defined on the probability space
$(\R^{\Z^{d}}, \B^{\Z^{d}}, \mu)$ is stationary with the same law as $(\Delta_{i})_{i\in\Z^{d}}$, hence, without loss of
generality, one can suppose that $(\Omega, \F, \P)=(\R^{\Z^{d}}, \B^{\Z^{d}}, \mu)$ and $\Delta_i=\pi_0\circ T^{i}$. Recall that $\rho$ is 
the metric defined for any finite subsets $B_1$ and $B_2$ of $\Z^d$ by $\rho(B_1,2_2)=\min\{\vert i-j\vert\,;\,i\in B_1,\,j\in B_2\}$ 
and $\vert i-j\vert=\max_{1\leq k\leq d}\vert i_k-j_k\vert$ for any $i=(i_1,...,i_d)$ and $j=(j_1,...,j_d)$ in $\Z^d$. We consider the following sets:
$$
\Lambda_{n}^{m_{n,\infty}}=\{i\in\Lambda_{n}\,;\,\rho(\{i\},\partial\Lambda_{n})\geq m_{n,\infty}\}
\quad\textrm{and}\quad \textrm{I}_{n}^{m_{n,\infty}}=\{1\leq k\leq \vert\Lambda_n\vert  \,;\,\varphi(k)\in\Lambda_{n}^{m_{n,\infty}}\},
$$
and the function $\Psi$ from $\R^{\Z^d}$ to $\R$ such that
$$
\Psi(\Delta)=\Delta_{0}^2+\sum_{i\in V_{0}^1\cap W_n}2\Delta_{0}\Delta_{i}\quad\textrm{where $W_n=\{-m_{n,\infty}+1,...,m_{n,\infty}-1\}^d$}.
$$
For $1\leq k\leq \vert\Lambda_n\vert  $, we set $D_{k}^{(n)}=\eta-b_n^{2\beta+1}\left(\Psi\circ T^{\varphi(k)}(\Delta)\right)$. By definition of $\Psi$ and of the set
$\textrm{I}_{n}^{m_{n,\infty}}$, we have for any $k$ in $\textrm{I}_{n}^{m_{n,\infty}}$
$$
\Psi\circ T^{\varphi(k)}(\Delta)=\Delta_{\varphi(k)}^2+2\Delta_{\varphi(k)}(S_{\varphi(k-1)}(\Delta)-S_{\varphi(k)}^{m_{n,\infty}}(\Delta)).
$$
Therefore for $k$ in $\textrm{I}_{n}^{m_{n,\infty}}$
$$
\frac{D_{k}^{(n)}}{\vert\Lambda_n\vert  }=\frac{\eta}{\vert\Lambda_n\vert  }-\Gamma_{\varphi(k)}^2-2\Gamma_{\varphi(k)}(S_{\varphi(k-1)}(\Gamma)-S_{\varphi(k)}^{m_{n,\infty}}(\Gamma)).
$$
Using ($\ref{condition_on_Lambda_n}$), we know that $\lim_{n\to+\infty}\vert\Lambda_n\vert^{-1}\vert\textrm{I}_{n}^{m_{n,\infty}}\vert=1$. So, it remains to consider
$$
\textrm{F}_1=\left\vert\E\left(\frac{1}{\vert\Lambda_n\vert  }\sum_{k=1}^{\vert\Lambda_n\vert  }h_{k-1,k+1}^{''}D_{k}^{(n)}\right)\right\vert.
$$
Applying Lemma $\ref{lemme-technique}$ and Lemma \ref{mn-tau}, we obtain
\begin{align*}
\textrm{F}_1&\leq\left\vert\E\left(\frac{b_n^{2\beta+1}}{\vert\Lambda_n\vert  }\sum_{k=1}^{\vert\Lambda_n\vert  }h_{k-1,k+1}^{''}(\Delta_{\varphi(k)}^2-\E(\Delta_0^2))\right)\right\vert
+\vert\eta-b_n^{2\beta+1}\E(\Delta_0^2)\vert+2b_n^{2\beta+1}\sum_{j\in V_0^1\cap W_n}\E\vert\Delta_0\Delta_{j}\vert\\
&\leq\left\vert\E\left(\frac{b_n^{2\beta+1}}{\vert\Lambda_n\vert}\sum_{k=1}^{\vert\Lambda_n\vert}h_{k-1,k+1}^{''}(\Delta_{\varphi(k)}^2-\E(\Delta_0^2))\right)\right\vert
+o(1).
\end{align*}
So, it suffices to prove that
$$
F_2=\left\vert\E\left(\frac{b_n^{2\beta+1}}{\vert\Lambda_n\vert  }\sum_{k=1}^{\vert\Lambda_n\vert  }h_{k-1,k+1}^{''}(\Delta_{\varphi(k)}^2-\E(\Delta_0^2))\right)\right\vert
$$
goes to zero as $n$ goes to infinity. Let $C>0$ be fixed. We have $F_2\leq F_2^{'}+F_2^{''}$ where
$$
F_2^{'}=\left\vert\E\left(\frac{b_n^{2\beta+1}}{\vert\Lambda_n\vert  }\sum_{k=1}^{\vert\Lambda_n\vert  }h_{k-1,k+1}^{''}
\left(\Delta_{\varphi(k)}^2-\E_C\left(\Delta_{\varphi(k)}^2\right)\right)\right)\right\vert
$$
and
$$
F_2^{''}=\left\vert\E\left(\frac{b_n^{2\beta+1}}{\vert\Lambda_n\vert  }\sum_{k=1}^{\vert\Lambda_n\vert  }h_{k-1,k+1}^{''}\left(\E_C\left(\Delta_{\varphi(k)}^2\right)-\E(\Delta_0^2)\right)\right)\right\vert
$$
where we recall the notation $\E_C\left(\Delta_{\varphi(k)}^2\right)=\E\left(\Delta_{\varphi(k)}^2\vert\F_{V_{\varphi(k)}^C}\right)$. 
The following result is a Serfling type inequality which can be found in \cite{McLeish}.
\begin{lemma}\label{Serfling-inequality}
Let $\U$ and $\mathcal{V}$ be two $\sigma$-algebras and let $X$ be a random variable measurable with respect to $\mathcal{U}$. 
If $1\leq p\leq r\leq\infty$ then 
$$
\|\E(X\vert\mathcal{V})-\E(X)\|_p\leq2(2^{1/p}+1)\left(\alpha(\U,\mathcal{V})\right)^{\frac{1}{p}-\frac{1}{r}}\|X\|_r.
$$
\end{lemma}
Applying Lemma $\ref{Serfling-inequality}$ and keeping in mind that $\vert\Delta_0\vert$ is bounded by $\kappa b_n^{-\beta-1}$, we derive 
$$
F_2^{''}\leq b_n^{2\beta+1}\|\E_C\left(\Delta_0^2\right)-\E(\Delta_0^2)\|_1\leq \kappa b_n^{-1}\alpha_{1,\infty}(C)
$$
In the other part,
$$
F_2^{'}\leq\frac{b_n^{2\beta+1}}{\vert\Lambda_n\vert  }\sum_{k=1}^{\vert\Lambda_n\vert  }\left(\textrm{J}_k^{1}(C)+\textrm{J}_k^2(C)\right)\\
$$
where
$$
\textrm{J}_k^1(C)=\left\vert\E\left(h_{k-1,k+1}^{''C}\circ T^{-\varphi(k)}
\left(\Delta_0^2-\E_C\left(\Delta_0^2\right)\right)\right)\right\vert=0
$$
since $h_{k-1,k+1}^{''C}\circ T^{-\varphi(k)}$ is $\F_{V_0^C}$-measurable and
\begin{align*}
b_n^{2\beta+1}\textrm{J}_k^2(C)&=b_n^{2\beta+1}\left\vert\E\left(\left(h_{k-1,k+1}^{''}\circ T^{-\varphi(k)}-h_{k-1,k+1}^{''C}\circ T^{-\varphi(k)}\right)
\left(\Delta_0^2-\E_C\left(\Delta_0^2\right)\right)\right)\right\vert\\
&\leq b_n^{2\beta+1}\E\left(\left(2\wedge\sum_{\vert i\vert<C}\frac{b_n^{\beta+\frac{1}{2}}\vert\Delta_i\vert}{\vert\Lambda_n\vert^{1/2}}\right)\Delta_0^2\right)\\
&\leq\frac{\kappa b_n^{2\beta+1}\E(\Delta_0^2)}{(\vert\Lambda_n\vert b_n)^{1/2}}+\frac{\kappa b_n^{2\beta+\frac{1}{2}}}{\vert\Lambda_n\vert^{1/2}}
\sum_{\substack{\vert i\vert<C \\ i\neq 0}}\E\vert\Delta_0\Delta_i\vert\quad\textrm{since $\vert\Delta_0\vert\leq \kappa b_n^{-\beta-1}$ a.s.}\\
&\leq \kappa\left(\frac{1}{(\vert\Lambda_n\vert b_n)^{1/2}}+\frac{C^d\sqrt{b_n}}{\vert\Lambda_n\vert^{1/2}}\right)\qquad\textrm{(by Lemma $\ref{lemme-technique}$)}
\end{align*}
So, putting $C=b_n^{\frac{-1}{2d-1}}$ and keeping in mind that $\sum_{m\geq0}m^{2d-1}\,\alpha_{1,\infty}(m)<+\infty$, we derive
$$
F_2=O\left(C^{2d-1}\,\alpha_{1,\infty}(C)\right)+O\left(\frac{1+b_n^{\frac{d-1}{2d-1}}}{(\vert\Lambda_n\vert  b_n)^{1/2}}\right)=o(1).
$$
The proof of Theorem $\ref{convergence-loi}$ is complete.
\subsection{Proof of Theorem $\textbf{\ref{convergence-loi-nonmixing-case}}$} 
The proof follows the same lines as in the proof of Theorem $\ref{convergence-loi}$. We consider the sequence $(m_n)_{n\geq 1}$ defined by 
\begin{equation}\label{def_mn}
m_n=\max\left\{v_n,\left[\left(\frac{1}{b_n^{3}}\sum_{\vert i\vert>v_n}\vert i\vert^{\frac{5d}{2}}\,\delta_i\right)^{\frac{1}{3d}}\right]+1\right\}
\end{equation}
where $v_n=\big[b_n^{\frac{-1}{2d}}\big]$ and $[\,.\,]$ denotes the integer part function. As in the proof of Theorem $\ref{convergence-loi}$, the sequence 
$(m_n)_{n\geq 1}$ satisfies the following lemma which the proof is left to the reader.
\begin{lemma}\label{mn}
If $\sum_{i\in\Z^d}\vert i\vert^{\frac{5d}{2}}\delta_i<+\infty$ holds then  
$$
m_n\to\infty,\quad m_n^db_n\to0\quad\textrm{and}\quad\frac{1}{(m_n^db_n)^{3/2}}\sum_{\vert i\vert>m_n}\vert i\vert^{\frac{5d}{2}}\,\delta_i\to 0.
$$
\end{lemma}
For any $z$ in $\R$, we denote
\begin{equation}\label{definition_Gnzi_overline_Gnzi}
\textrm{G}_n(z,i)=g_n\left(\frac{z-Y_i}{b_n}\right)\quad\textrm{and}\quad
\overline{\textrm{G}}_n(z,i)=\E\left(\textrm{G}_n(z,i)\vert\F_{n,i}\right)
\end{equation}
where $\F_{n,i}=\sigma\left(\varepsilon_{i-s}\,;\,\vert s\vert\leq m_n\right)$. 
Denoting $M_n=2m_n+1$, $(\overline{\textrm{G}}_n(z,i))_{i\in\Z^d}$ is an $M_n$-dependent random field 
(i.e. $\overline{\textrm{G}}_n(z,i)$ and $\overline{\textrm{G}}_n(z,j)$ are independent as soon as $\vert i-j\vert\geq M_n$). 
For any $x$ in $\R$ and any integer $n\geq 1$, we denote
$$
\overline{f}_n(x)=\frac{1}{\vert\Lambda_n\vert b_n}\sum_{i\in\Lambda_n}\overline{\textrm{G}}_n(x,i).
$$ 
The proof of the following lemma is done in the appendix.
\begin{lemma}\label{moment_inequality}
For any $p\geq 2$, any $x$ in $\R$, any positive integer $n$ and any $(a_i)_{i\in\Z^d}$ in $\R^{\Z^d}$, 
$$
\left\|\sum_{i\in\Lambda_n}a_i\left(\emph{\textrm{G}}_n(x,i)-\overline{\emph{\textrm{G}}}_n(x,i)\right)\right\|_p
\leq\frac{\kappa m_n^d}{b_n^{1+\beta}}\left(p\sum_{i\in\Lambda_n}a_i^2\right)^{1/2}\sum_{\vert i\vert>m_n}\delta_{i,p}.
$$
\end{lemma}
Let $x\neq y$ be fixed and let $\lambda_1$ and $\lambda_2$ be two constants such that $\lambda_1^2+\lambda_2^2=1$. We have  
\begin{align*}
\lambda_1(\vert\Lambda_n\vert b_n^{2\beta+1})^{1/2}(\hat{f}_n(x)-\E\hat{f}_n(x))+\lambda_2(\vert\Lambda_n\vert b_n^{2\beta+1})^{1/2}(\hat{f}_n(y)-\E\hat{f}_n(y))
&=\sum_{i\in\Lambda_n}\frac{b_n^{\beta+\frac{1}{2}}\Delta_i}{\vert\Lambda_n\vert^{1/2}}\\
\lambda_1(\vert\Lambda_n\vert b_n^{2\beta+1})^{1/2}(\overline{f}_n(x)-\E\overline{f}_n(x))+\lambda_2(\vert\Lambda_n\vert b_n^{2\beta+1})^{1/2}(\overline{f}_n(y)-\E\overline{f}_n(y))
&=\sum_{i\in\Lambda_n}\frac{b_n^{\beta+\frac{1}{2}}\overline{\Delta}_i}{\vert\Lambda_n\vert^{1/2}}
\end{align*} 
where 
$$
\Delta_i=\lambda_1 Z_{i}(x)+\lambda_2 Z_{i}(y)\quad\textrm{and}\quad\overline{\Delta}_i=\lambda_1\overline{Z}_{i}(x)+\lambda_2\overline{Z}_{i}(y)
$$
and for any $z$ in $\R$,
$$
Z_{i}(z)=\frac{1}{b_n}\left(\textrm{G}_n(z,i)-\E\textrm{G}_n(z,i)\right)
\quad\textrm{and}\quad\overline{Z}_{i}(z)=\frac{1}{b_n}\left(\overline{\textrm{G}}_n(z,i)-\E\overline{\textrm{G}}_n(z,i)\right)
$$
where $\textrm{G}_n(z,i)$ and $\overline{\textrm{G}}_n(z,i)$ are defined by $(\ref{definition_Gnzi_overline_Gnzi})$. 
Applying Lemma \ref{mn} and Lemma \ref{moment_inequality}, we know that
\begin{equation}\label{ecart_sommes_partielles_delta_et_overline_delta}
\frac{b_n^{\beta+\frac{1}{2}}}{\vert\Lambda_n\vert^{1/2}}\left\|\sum_{i\in\Lambda_n}\left(\Delta_i-\overline{\Delta}_i\right)\right\|_2
\leq\frac{\kappa (\vert\lambda_1\vert+\vert\lambda_2\vert)}{(m_n^db_n)^{3/2}}\sum_{\vert i\vert>m_n}\vert i\vert^{\frac{5d}{2}}\delta_i=o(1).
\end{equation}
So, it suffices to prove the asymptotic normality of $\left(b_n^{\beta+\frac{1}{2}}\vert\Lambda_n\vert^{-1/2}\sum_{i\in\Lambda_n}\overline{\Delta}_i\right)_{n\geq 1}$. Let $\eta$ be defined by $(\ref{eta_sigma})$. 
The proof of the following lemma is also postponed to the appendix.
\begin{lemma}\label{lemme-technique-barre} $b_n^{2\beta+1}\E(\overline{\Delta}_0^2)$ converges to $\eta$ as $n$ goes to infinity. Moreover, if ($\ref{nonmixing-condition}$) holds then
$\sup_{i\in\Z^d\backslash\{0\}}\E\vert\overline{\Delta}_0\overline{\Delta}_i\vert=o(m_n^{-d}b_n^{-2\beta-1})$.
\end{lemma}
As in the proof of Theorem \ref{convergence-loi}, in order to describe the set $\Lambda_n$, we consider the one to one map 
$\varphi$ from $[1,\vert\Lambda_n\vert]\cap\N^{\ast}$ to $\Lambda_{n}$ by:
$\varphi$ is the unique function such that $\varphi(k)<_{\text{\text{lex}}}\varphi(l)$ for $1\leq k<l\leq
\vert\Lambda_n\vert$ where $<_{\textrm{lex}}$ denotes the lexicographic order on $\Z^d$ and we consider a field 
$(\tau_{i})_{i\in\Z^d}$ of i.i.d. random variables independent of $(\overline{\Delta}_{i})_{i\in\Z^d}$ such that $\tau_{0}$ has the standard normal
law $\mathcal{N}(0,1)$. We introduce the fields $\Gamma$ and $\gamma$ defined for any $i$ in $\Z^d$ by 
$$
\Gamma_{i}=\frac{b_n^{\beta+\frac{1}{2}}\overline{\Delta}_i}{\vert\Lambda_n\vert^{1/2}}\quad\textrm{and}\quad\gamma_{i}=\frac{\tau_i\sqrt{\eta}}{\vert\Lambda_n\vert^{1/2}}.
$$ 
Note that $\Gamma$ is an $M_n$-dependent random field where $M_n=2m_n+1$ and $m_n$ is defined by (\ref{def_mn}). 
Keeping in mind the notations introduced in the proof of Theorem \ref{convergence-loi}, it suffices to prove that for any function $h$ in 
$B_{1}^4(\R)$,
$$
\E\left(h\left(S_{\varphi(\vert\Lambda_n\vert)}(\Gamma)\right)\right)\converge{n}{+\infty}{}\E \left(h\left(\tau_0\sqrt{\eta}\right)\right).
$$
Applying Lindeberg's decomposition, we have
$$
\E \left(h\left(S_{\varphi(\vert\Lambda_n\vert)}(\Gamma)\right)-h\left(\tau_{0}\sqrt{\eta}\right)\right)
=\sum_{k=1}^{\vert\Lambda_n\vert}\E \left(h_{k,k+1}-h_{k-1,k}\right).
$$
Now,
$$
h_{k,k+1}-h_{k-1,k}=h_{k,k+1}-h_{k-1,k+1}+h_{k-1,k+1}-h_{k-1,k}.
$$
By Taylor's formula we get
$$
h_{k,k+1}-h_{k-1,k+1}=\Gamma_{\varphi(k)}h_{k-1,k+1}^{'}+\frac{1}{2}\Gamma_{\varphi(k)}^{2}h_{k-1,k+1}^{''}+R_{k}
$$
and
$$
h_{k-1,k+1}-h_{k-1,k}=-\gamma_{\varphi(k)}h_{k-1,k+1}^{'}-\frac{1}{2}\gamma_{\varphi(k)}^{2}h_{k-1,k+1}^{''}+r_{k}
$$
where $\vert R_{k}\vert\leq \Gamma_{\varphi(k)}^2(1\wedge\vert \Gamma_{\varphi(k)}\vert)$ and 
$\vert r_{k}\vert\leq\gamma_{\varphi(k)}^2(1\wedge\vert\gamma_{\varphi(k)}\vert)$.
Since $(\Gamma,\tau_{i})_{i\neq \varphi(k)}$ is independent of $\tau_{\varphi(k)}$, it follows that
$$
\E \left(\gamma_{\varphi(k)}h_{k-1,k+1}^{'}\right)=0\quad\textrm{and}\quad
\E \left(\gamma_{\varphi(k)}^2h_{k-1,k+1}^{''}\right)=\E \left(\frac{\eta}{\vert\Lambda_n\vert}h_{k-1,k+1}^{''}\right)
$$
Hence, we obtain
\begin{align*}
\E \left(h(S_{\varphi(\vert\Lambda_n\vert)}(\Gamma))-h\left(\tau_0\sqrt{\eta}\right)\right)&=
\sum_{k=1}^{\vert\Lambda_n\vert}\E (\Gamma_{\varphi(k)}h_{k-1,k+1}^{'})\\
&\quad+\sum_{k=1}^{\vert\Lambda_n\vert}\E \left(\left(\Gamma_{\varphi(k)}^2-\frac{\eta}{\vert\Lambda_n\vert}\right)\frac{h_{k-1,k+1}^{''}}{2}\right)\\
&\quad+\sum_{k=1}^{\vert\Lambda_n\vert}\E \left(R_{k}+r_{k}\right).
\end{align*}
Let $1\leq k\leq \vert\Lambda_n\vert$ be fixed. Noting that $\vert\overline{\Delta}_0\vert$ is bounded by $\kappa b_n^{-\beta-1}$ and applying Lemma $\ref{lemme-technique-barre}$, 
we derive
$$
\E\vert R_k\vert\leq\frac{b_n^{3\beta+\frac{3}{2}}\E\vert\overline{\Delta}_0\vert^3}{\vert\Lambda_n\vert^{3/2}}=O\left(\frac{1}{(\vert\Lambda_n\vert^{3}\,b_n)^{1/2}}\right)
$$
and
$$
\E\vert r_k\vert\leq\frac{\E\vert\gamma_0\vert^3}{\vert\Lambda_n\vert^{3/2}}\leq\frac{\eta^{3/2}\E\vert\tau_0\vert^3}{\vert\Lambda_n\vert^{3/2}}
=O\left(\frac{1}{\vert\Lambda_n\vert^{3/2}}\right).
$$
Consequently, we obtain
$$
\sum_{k=1}^{\vert\Lambda_n\vert}\E \left(\vert R_{k}\vert+\vert r_{k}\vert\right)=O\left(\frac{1}{(\vert\Lambda_n\vert b_n)^{1/2}}+\frac{1}{\vert\Lambda_n\vert^{1/2}}\right)
=o(1).
$$
Now, it is sufficient to show
\begin{equation}\label{equation1_nonmixing_case}
\lim_{n\to+\infty}\sum_{k=1}^{\vert\Lambda_n\vert}\left(\E (\Gamma_{\varphi(k)}h_{k-1,k+1}^{'})+\E \left(\left(\Gamma_{\varphi(k)}^2-\frac{\eta}{\vert\Lambda_n\vert}\right)\frac{h_{k-1,k+1}^{''}}{2}\right)\right)=0.
\end{equation}
We focus on $\sum_{k=1}^{\vert\Lambda_n\vert}\E \left(\Gamma_{\varphi(k)}h_{k-1,k+1}^{'}\right)$. 
Recall that the sets $\{V_{i}^{k}\,;\,i\in\Z^{d}\,,\,k\in\N^{\ast}\}$ are defined as follows: 
$$
V_{i}^{1}=\{j\in\Z^{d}\,;\,j<_{\textrm{lex}}i\}\,\,\,\textrm{and for $k\geq 2$,}\,\, V_{i}^{k}=V_{i}^{1}\cap\{j\in\Z^{d}\,;\,\vert i-j\vert\geq k\}. 
$$
For all $n$ in $\N^{\ast}$ and all integer $k$ in $[1,\vert\Lambda_n\vert]$, we define
$$
\textrm{E}_{k}^{M_n}=\varphi([1,k]\cap\N^{\ast})\cap V_{\varphi(k)}^{M_n}\quad\textrm{and}\quad
S_{\varphi(k)}^{M_n}(\Gamma)=\sum_{i\in\textrm{E}_{k}^{M_n}}\Gamma_{i}
$$
where $M_n=2m_n+1$. For any function $\Psi$ from $\R$ to $\R$, we define $\Psi_{k-1,l}^{M_n}=\Psi(S_{\varphi(k)}^{M_n}(\Gamma)+S_{\varphi(l)}^c(\gamma))$. 
Our aim is to show that
\begin{equation}\label{equation1bis}
\lim_{n\to+\infty}\sum_{k=1}^{\vert\Lambda_n\vert}\E\left(\Gamma_{\varphi(k)}h_{k-1,k+1}^{'}-\Gamma_{\varphi(k)}
\left(S_{\varphi(k-1)}(\Gamma)-S_{\varphi(k)}^{M_n}(\Gamma)\right)h_{k-1,k+1}^{''}\right)=0.
\end{equation}
First, we note that
$$
\Gamma_{\varphi(k)}h_{k-1,k+1}^{'}=\Gamma_{\varphi(k)}h_{k-1,k+1}^{'M_n}+\Gamma_{\varphi(k)}\left(h_{k-1,k+1}^{'}-h_{k-1,k+1}^{'M_n}\right).
$$
Applying again Taylor's formula,
$$
\Gamma_{\varphi(k)}(h_{k-1,k+1}^{'}-h_{k-1,k+1}^{'M_n})=\Gamma_{\varphi(k)}(S_{\varphi(k-1)}(\Gamma)-S_{\varphi(k)}^{M_n}(\Gamma))h_{k-1,k+1}^{''}+R_{k}^{'},
$$
where 
$$
\vert R_{k}^{'}\vert\leq 2\vert \Gamma_{\varphi(k)}(S_{\varphi(k-1)}(\Gamma)-S_{\varphi(k)}^{M_n}(\Gamma))(1\wedge\vert S_{\varphi(k-1)}(\Gamma)-S_{\varphi(k)}^{M_n}(\Gamma)\vert)\vert.
$$ 
Since $(\Gamma_i)_{i\in\Z^d}$ is $M_n$-dependent, we have $\E \left(\Gamma_{\varphi(k)}h_{k-1,k+1}^{'{M_n}}\right)=0$ 
and consequently $(\ref{equation1bis})$ holds if and only if $\lim_{n\to+\infty}\sum_{k=1}^{\vert\Lambda_n\vert}\E \vert R_{k}^{'}\vert=0$. In fact, 
denoting $W_n=\{-M_n+1,...,M_n-1\}^d$ and $W_n^{\ast}=W_n\backslash\{0\}$, it follows that
\begin{align*}
\sum_{k=1}^{\vert\Lambda_n\vert}\E \vert R_{k}^{'}\vert 
&\leq 2b_n^{2\beta+1}\E \left(\vert\overline{\Delta}_{0}\vert\left(\sum_{i\in W_n}\vert\overline{\Delta}_{i}\vert\right)
\left(1\wedge\frac{b_n^{\beta+\frac{1}{2}}}{\vert\Lambda_n\vert^{1/2}}\sum_{i\in W_n}\vert\overline{\Delta}_{i}\vert\right)\right)\\
&=2b_n^{2\beta+1}\E\left(\left(\overline{\Delta}_0^2+\sum_{i\in W_n^{\ast}}\vert\overline{\Delta}_0\overline{\Delta}_i\vert\right)\left(1\wedge\frac{b_n^{\beta+\frac{1}{2}}}{\vert\Lambda_n\vert^{1/2}}\sum_{i\in W_n}\vert\overline{\Delta}_i\vert\right)\right)\\
&\leq\frac{2b_n^{3\beta+\frac{3}{2}}}{\vert\Lambda_n\vert^{1/2}}\sum_{i\in W_n}\E(\overline{\Delta}_0^2\vert\overline{\Delta}_i\vert)+2b_n^{2\beta+1}\sum_{i\in W_n^{\ast}}\E\vert\overline{\Delta}_0\overline{\Delta}_i\vert.
\end{align*}
Since $\vert\overline{\Delta}_0\vert$ is bounded by $\kappa b_n^{-\beta-1}$, we derive
\begin{align*}
\sum_{k=1}^{\vert\Lambda_n\vert}\E \vert R_{k}^{'}\vert 
&\leq\frac{\kappa b_n^{2\beta+\frac{1}{2}}}{\vert\Lambda_n\vert^{1/2}}\sum_{i\in W_n}\E(\vert\overline{\Delta}_0\overline{\Delta}_i\vert)+2b_n^{2\beta+1}\sum_{i\in W_n^{\ast}}\E\vert\overline{\Delta}_0\overline{\Delta}_i\vert\\
&\leq\frac{\kappa b_n^{2\beta+1}\E (\overline{\Delta}_0^2)}{(\vert\Lambda_n\vert b_n)^{1/2}}+\kappa b_n^{2\beta+1}\left(1+\frac{1}{(\vert\Lambda_n\vert b_n)^{1/2}}\right)m_n^d\sup_{i\in\Z^d\backslash\{0\}}\E(\vert\overline{\Delta}_0\overline{\Delta}_i\vert)\\
&=o(1)\quad\textrm{(by Lemma $\ref{lemme-technique-barre}$ and Assumption $\textbf{(A5)}$}).
\end{align*}
In order to obtain ($\ref{equation1_nonmixing_case}$) it remains to control
$$
F_{0}=\E \left(\sum_{k=1}^{\vert\Lambda_n\vert}h_{k-1,k+1}^{''}
\left(\frac{\Gamma_{\varphi(k)}^2}{2}+\Gamma_{\varphi(k)}\left(S_{\varphi(k-1)}(\Gamma)-S_{\varphi(k)}^{M_n}(\Gamma)\right)-\frac{\eta}{2\vert\Lambda_n\vert}\right)\right).
$$
Denote by $\B$ the Borel $\sigma$-algebra of $\R$. Without loss of generality, one can suppose that $(\Omega, \F, \P)=(\R^{\Z^{d}}, \B^{\Z^{d}}, \mu)$ 
and $\overline{\Delta}_{k}=\pi_0\circ T^{k}$ where 
$\mu$ is the law of the stationary real random field $(\overline{\Delta}_{k})_{k\in\Z^{d}}$, $\pi_0$ is the projection from $\R^{\Z^{d}}$ to $\R$ 
defined by $\pi_0(\omega)=\omega_{0}$ and $(T^{k})_{k\in\Z^{d}}$ is the family of translation operators from $\R^{\Z^{d}}$ to $\R^{\Z^{d}}$ defined by
$(T^{k}(\omega))_{i}=\omega_{i+k}$ for any $k\in\Z^{d}$ and any $\omega$ in $\R^{\Z^{d}}$. Recall also the metric 
$\rho$ defined for any finite subsets $B_1$ and $B_2$ of $\Z^d$ by $\rho(B_1,2_2)=\min\{\vert i-j\vert\,;\,i\in B_1,\,j\in B_2\}$ 
and $\vert i-j\vert=\max_{1\leq k\leq d}\vert i_k-j_k\vert$ for any $i=(i_1,...,i_d)$ and $j=(j_1,...,j_d)$ in $\Z^d$. We consider the following sets:
$$
\Lambda_{n}^{M_n}=\{i\in\Lambda_{n}\,;\,\rho(\{i\},\partial\Lambda_{n})\geq M_n\}
\quad\textrm{and}\quad \I_{n}^{M_n}=\{1\leq k\leq \vert\Lambda_n\vert\,;\,\varphi(k)\in\Lambda_{n}^{M_n}\},
$$
and the function $\Psi$ from $\R^{\Z^d}$ to $\R$ such that
$$
\Psi(\overline{\Delta})=\overline{\Delta}_{0}^2+\sum_{i\in V_{0}^1\cap W_n}2\overline{\Delta}_{0}\overline{\Delta}_{i}.
$$
For $1\leq k\leq \vert\Lambda_n\vert$, we set $D_{k}^{(n)}=\eta-b_n^{2\beta+1}\Psi\circ T^{\varphi(k)}(\overline{\Delta})$. By definition of $\Psi$ and of the set
$\I_{n}^{M_n}$, we have for any $k$ in $\I_{n}^{M_n}$, 
$\Psi\circ T^{\varphi(k)}(\overline{\Delta})=\overline{\Delta}_{\varphi(k)}^2
+2\overline{\Delta}_{\varphi(k)}(S_{\varphi(k-1)}(\overline{\Delta})-S_{\varphi(k)}^{M_n}(\overline{\Delta}))$.
Therefore for $k$ in $\I_{n}^{M_n}$,
$$
\frac{D_{k}^{(n)}}{\vert\Lambda_n\vert}=\frac{\eta}{\vert\Lambda_n\vert}-\Gamma_{\varphi(k)}^2-2\Gamma_{\varphi(k)}(S_{\varphi(k-1)}(\Gamma)-S_{\varphi(k)}^{M_n}(\Gamma)).
$$
Since $(\ref{condition_on_Lambda_n})$ holds, we have $\lim_{n\to+\infty}\vert\Lambda_n\vert^{-1}\vert \I_{n}^{M_n}\vert=1$. So, it remains to consider
$$
F_1=\left\vert\E\left(\frac{1}{\vert\Lambda_n\vert}\sum_{k=1}^{\vert\Lambda_n\vert}h_{k-1,k+1}^{''}D_{k}^{(n)}\right)\right\vert.
$$
Applying Lemma $\ref{lemme-technique-barre}$, we have
\begin{align*}
F_1&\leq\left\vert\E\left(\frac{b_n^{2\beta+1}}{\vert\Lambda_n\vert}\sum_{k=1}^{\vert\Lambda_n\vert}h_{k-1,k+1}^{''}(\overline{\Delta}_{\varphi(k)}^2-\E(\overline{\Delta}_0^2))\right)\right\vert
+\vert\eta-b_n^{2\beta+1}\E(\overline{\Delta}_0^2)\vert+2b_n^{2\beta+1}\sum_{j\in V_0^1\cap W_n}\E\vert\overline{\Delta}_0\overline{\Delta}_{j}\vert\\
&\leq\left\vert\E\left(\frac{b_n^{2\beta+1}}{\vert\Lambda_n\vert}\sum_{k=1}^{\vert\Lambda_n\vert}h_{k-1,k+1}^{''}(\overline{\Delta}_{\varphi(k)}^2-\E(\overline{\Delta}_0^2))\right)\right\vert
+o(1).
\end{align*}
So, it suffices to prove that
$$
F_2=\left\vert\E\left(\frac{b_n^{2\beta+1}}{\vert\Lambda_n\vert}\sum_{k=1}^{\vert\Lambda_n\vert}h_{k-1,k+1}^{''}(\overline{\Delta}_{\varphi(k)}^2-\E(\overline{\Delta}_0^2))\right)\right\vert
$$
goes to zero as $n$ goes to infinity. In fact, 
$F_2\leq\frac{b_n^{2\beta+1}}{\vert\Lambda_n\vert}\sum_{k=1}^{\vert\Lambda_n\vert}\left(\textrm{J}_k^{1}(n)+\textrm{J}_k^2(n)\right)$ where 
$\textrm{J}_k^1(n)=\left\vert\E\left(h_{k-1,k+1}^{''M_n}\left(\overline{\Delta}_{\varphi(k)}^2-\E\left(\overline{\Delta}_0^2\right)\right)\right)\right\vert=0$ 
since $h_{k-1,k+1}^{''M_n}$ is $\sigma\left(\overline{\Delta}_i\,;\,i\in V_{\varphi(k)}^{M_n}\right)$-measurable and the conditional expectation of 
$\overline{\Delta}_{\varphi(k)}^2$ with respect to $\sigma\left(\overline{\Delta}_i\,;\,i\in V_{\varphi(k)}^{M_n}\right)$ is null 
(cf. $(\overline{\Delta}_i)_{i\in\Z^d}$ is an $M_n$-dependent random field).
\begin{align*}
b_n^{2\beta+1}\textrm{J}_k^2(n)
&=b_n^{2\beta+1}\left\vert\E\left(\left(h_{k-1,k+1}^{''}-h_{k-1,k+1}^{''M_n}\right)
\left(\overline{\Delta}_{\varphi(k)}^2-\E\left(\overline{\Delta}_0^2\right)\right)\right)\right\vert\\
&\leq b_n^{2\beta+1}\E\left(\left(2\wedge\sum_{\vert i\vert<M_n}\frac{b_n^{\beta+\frac{1}{2}}\vert\overline{\Delta}_i\vert}{\vert\Lambda_n\vert^{1/2}}\right)\overline{\Delta}_0^2\right)\\
&\leq\frac{\kappa b_n^{2\beta+1}\E(\overline{\Delta}_0^2)}{(\vert\Lambda_n\vert b_n)^{1/2}}+\frac{\kappa b_n^{2\beta+\frac{1}{2}}}{\vert\Lambda_n\vert^{1/2}}
\sum_{\substack{\vert i\vert<M_n \\ i\neq 0}}\E\vert\overline{\Delta}_0\overline{\Delta}_i\vert\quad\textrm{since $\vert\overline{\Delta}_0\vert\leq \kappa b_n^{-\beta-1}$ a.s.}\\
&\leq\frac{\kappa\left(b_n^{2\beta+1}\E(\overline{\Delta}_0^2)+b_n^{2\beta+1}m_n^d\sup_{i\in\Z^d\backslash\{0\}}\E(\vert\overline{\Delta}_0\overline{\Delta}_i\vert)\right)}{(\vert\Lambda_n\vert b_n)^{1/2}}\\
&=o(1)\quad\textrm{(by Lemma $\ref{lemme-technique-barre}$ and Assumption $\textbf{(A5)}$)}.
\end{align*}
The proof of Theorem $\ref{convergence-loi-nonmixing-case}$ is complete.
\section{Appendix}
{\em Proof of Lemma $\ref{mn-tau}$}. Let $\tau$ be fixed in $\N^{\ast}\cup\{\infty\}$ and let $m_{n,\tau}$ be defined by equation ($\ref{mn-tau}$). 
Since $b_n$ goes to zero and $m_{n,\tau}\geq v_n=[b_n^{-\frac{1}{2d}}]$ where $[\,.\,]$ is the integer part function, we obtain $m_{n,\tau}$ goes to infinity as 
$n$ goes to infinity. For any positive integer $m$, we consider
$$
r(m)=\sum_{\vert i\vert>m}\vert i\vert^d\,\alpha_{1,\tau}(\vert i\vert).
$$
Since $\sum_{m\geq 1}m^{2d-1}\,\alpha_{1,\tau}(m)<\infty$, we have $r(m)$ converges to zero as $m$ goes to infinity. Moreover,
$$
m_{n,\tau}^db_n\leq\max\left\{\sqrt{b_n},\kappa\left(\sqrt{r\left(v_n\right)}+b_n\right)\right\}\converge{n}{+\infty}{ }0
$$
We have also
$$
m_{n,\tau}^d\geq\frac{1}{b_n}\sqrt{r\left(v_n\right)}\geq\frac{1}{b_n}\sqrt{r\left(m_{n,\tau}\right)}\qquad\textrm{since $v_n\leq m_{n,\tau}$}.
$$
Finally, we obtain
$$
\frac{1}{m_{n,\tau}^db_n}\sum_{\vert i\vert>m_{n,\tau}}\vert i\vert^d\,\alpha_{1,\tau}(\vert i\vert)\leq\sqrt{r(m_{n,\tau})}\converge{n}{+\infty}{ }0.
$$
The proof of Lemma $\ref{mn-tau}$ is complete.\\
\\
{\em Proof of Lemma $\ref{lemme_Z}$}. Using Assumption $\textbf{(A3)}$, for any $t$ in $\R$, there exists a sequence $(r_n(t))_{n\geq 1}$ 
going to $1$ as $n$ goes to infinity such that $\phi_{\theta}(t/b_n)=Br_n(t)b_n^{\beta}/\vert t\vert^{\beta}$. Consequently,
\begin{align*}
\int_{\R}g_n^2(u)du
&=\int_{\R}\left(\frac{1}{2\pi}\int_{\R}e^{-itu}\frac{\phi_K(t)}{\phi_{\theta}(t/b_n)}dt\right)^2du\\
&=\int_{\R}\left(\frac{1}{2\pi}\int_{\R}e^{-itu}\frac{\vert t\vert^{\beta}\phi_K(t)}{Br_n(t)b_n^{\beta}}dt\right)^2du\\
&=\frac{1}{B^2b_n^{2\beta}}\int_{\R}\frac{\vert t\vert^{2\beta}\phi_K^2(t)}{r_n^2(t)}dt\quad\textrm{(By Plancherel's theorem)}\\
\end{align*}
So, we obtain 
\begin{equation}\label{limite_integrale_carre_gn}
\lim_{n\to\infty}b_n^{2\beta}\int_{\R}g_n^2(u)du=\frac{1}{B^2}\int_{\R}\vert t\vert^{2\beta}\phi_K^2(t)dt
\end{equation}
Similarly, 
\begin{equation}\label{limite_integrale_abs_gn}
\lim_{n\to\infty}b_n^{\beta}\int_{\R}\vert g_n(u)\vert du=\frac{1}{B}\int_{\R}\vert t\vert^{\beta}\vert\phi_K(t)\vert dt
\end{equation}
For any $z$ in $\R$ and any $i$ in $\Z^d$, we recall that $\textrm{G}_n(z,i):=g_n\left(\frac{z-Y_i}{b_n}\right)$. 
Using ($\ref{limite_integrale_abs_gn}$) and keeping in mind Assumption $\textbf{(A1)}$, we have 
\begin{equation}\label{esperance_gn}
b_n^{\beta-1}\E\left\vert\textrm{G}_n(z,i)\right\vert=b_n^{\beta}\int_{\R} \vert g_n(u)\vert f_Y(z-ub_n)du\converge{n}{+\infty}{ }
\frac{f_Y(z)}{B}\int_{\R}\vert t\vert^{\beta}\vert\phi_K(t)\vert dt. 
\end{equation}
Similarly, using ($\ref{limite_integrale_carre_gn}$), we derive
\begin{equation}\label{variance_gn}
b_n^{2\beta-1}\E\textrm{G}^2_n(z,i)=b_n^{2\beta}\int_{\R} g_n^2(u)f_Y(z-ub_n)du\converge{n}{+\infty}{ }\frac{f_Y(z)}{B^2}\int_{\R}\vert t\vert^{2\beta}\phi_K^2(t)dt.
\end{equation}
In the other part, for any $s$ in $\R$,
\begin{equation}\label{decomposition_variance_Z}
b_n^{2\beta+1}\E(Z_0^2(s))=b_n^{2\beta-1}\left(\E\textrm{G}_n^2(s,0)-\left(\E\textrm{G}_n(s,0)\right)^2\right).
\end{equation}
Combining ($\ref{esperance_gn}$), ($\ref{variance_gn}$) and ($\ref{decomposition_variance_Z}$), we derive for any $s$ in $\R$,
\begin{equation}\label{convergence_esperance_Zs_carre_bis}
\lim_{n\to\infty}b_n^{2\beta+1}\E(Z_0^2(s))=\frac{f_Y(s)}{B^2}\int_{\R}\vert t\vert^{2\beta}\phi_K^2(t)dt.
\end{equation}
Now, let $s$ and $t$ be fixed in $\R$ and let $i$ be fixed in $\Z^d\backslash\{0\}$. Then,
\begin{equation}\label{inequality_esperance_Z0_Zi}
b_n^2\E\vert Z_0(s)Z_i(t)\vert
\leq\E\left\vert\textrm{G}_n(s,0)\textrm{G}_n(t,i)\right\vert
+3\E\left\vert \textrm{G}_n(s,0)\right\vert\E\left\vert \textrm{G}_n(t,0)\right\vert
\end{equation}
and
\begin{align*}
\E\left\vert \textrm{G}_n(s,0)\textrm{G}_n(t,i)\right\vert
&\leq b_n^2\iint \vert g_n(u)g_n(v)\vert f_Y(s-ub_n)\left\vert f_{i\vert 0}(s-ub_n,t-vb_n)-f_Y(t-vb_n)\right\vert dudv\\
&\qquad+b_n^2 \int\vert g_n(u)\vert f_Y(s-ub_n)du\int\vert g_n(v)\vert f_Y(t-vb_n)dv
\end{align*}
where $f_{i\vert 0}$ is the conditional density of $Y_i$ given $Y_0$. Using $(\ref{esperance_gn})$ and Assumption $\textbf{(A2)}$, we obtain 
\begin{equation}\label{covariance_de_gn}
\sup_{i\in\Z^d\backslash\{0\}}\E\left\vert \textrm{G}_n(s,0)\textrm{G}_n(t,i)\right\vert\leq \kappa b_n^{2-2\beta}.
\end{equation}
Combining $(\ref{esperance_gn})$, $(\ref{inequality_esperance_Z0_Zi})$ and $(\ref{covariance_de_gn})$, 
we derive 
$$
\sup_{i\in\Z^d\backslash\{0\}}\E\vert Z_0(s)Z_i(t)\vert=O(b_n^{-2\beta})
$$ 
The proof of Lemma \ref{lemme_Z} is complete.\\
\\
{\em Proof of Lemma $\ref{lemme-technique}$}. For any $z$ in $\R$ and any $i$ in $\Z^d$, recall that $\textrm{G}_n(z,i)=g_n\left(\frac{z-Y_i}{b_n}\right)$. We have
\begin{equation}\label{decomposition_esperance_delta_0_carre}
\E(\Delta_0^2)=\lambda_1^2\E(Z_0^2(x))+\lambda_2^2\E(Z_0^2(y))+2\lambda_1\lambda_2\E(Z_0(x)Z_0(y))
\end{equation}
and
\begin{equation}\label{decomposition_Zx_Zy}
b_n^2\E(Z_0(x)Z_0(y))=\E\left(\textrm{G}_n(x,0)\textrm{G}_n(y,0)\right)-\E\left(\textrm{G}_n(x,0)\right)\E\left(\textrm{G}_n(y,0)\right).
\end{equation}
Moreover,
$$
\E\left(\textrm{G}_n(x,0)\textrm{G}_n(y,0)\right)=b_n\int_{\R}g_n(u)g_n\left(u+\frac{y-x}{b_n}\right)f_Y(x-ub_n)du.
$$
Keeping in mind Assumption $\textbf{(A4)}$ and $x\neq y$, we have
$$
\left\vert g_n\left(u+\frac{y-x}{b_n}\right)\right\vert
\leq\frac{1}{2\pi\left\vert u+(y-x)/b_n\right\vert}\int_{\R}\left\vert\frac{\Phi_K\left(\frac{t}{u+(y-x)/b_n}\right)}{\Phi_{\theta}\left(\frac{t}{y-x+ub_n}\right)}\right\vert dt
=O(b_n).
$$
So, using ($\ref{esperance_gn})$, we derive 
\begin{equation}\label{esperance_gnx_gny}
\vert\E\left(\textrm{G}_n(x,0)\textrm{G}_n(y,0)\right)\vert\leq\kappa b_n^2\int_{\R}\vert g_n(u)\vert f_Y(x-ub_n)du=O(b_n^{2-\beta}).
\end{equation}
Combining ($\ref{esperance_gn}$), ($\ref{decomposition_Zx_Zy}$) and ($\ref{esperance_gnx_gny}$) and applying Lemma \ref{mn-tau}, we derive 
\begin{equation}\label{limite_esperance_Zx_Zy}
\vert\E(Z_0(x)Z_0(y))\vert=O(b_n^{-2\beta}).
\end{equation}
Combining ($\ref{convergence_esperance_Zs_carre_bis}$), $(\ref{decomposition_esperance_delta_0_carre})$ and $(\ref{limite_esperance_Zx_Zy}$), 
we obtain $\lim_{n\to+\infty}b_n^{2\beta+1}\E(\Delta_0^2)=\eta$. Let $i$ in $\Z^d\backslash\{0\}$ be fixed. Noting that 
$$
\E\vert\Delta_0\Delta_i\vert\leq\lambda_1^2\E\vert Z_0(x)Z_i(x)\vert+2\lambda_1\lambda_2\E\vert Z_0(x)Z_i(y)\vert+\lambda_2^2\E\vert Z_0(y)Z_i(y)\vert
$$
and applying the second part of Lemma \ref{lemme_Z}, we obtain $\sup_{i\in\Z^d\backslash\{0\}}\E\vert\Delta_0\Delta_i\vert=O(b_n^{-2\beta})$. 
The proof of Lemma \ref{lemme-technique} is complete.\\
\\
{\em Proof of Lemma \ref{moment_inequality}}. We follow the proof of Proposition 1 in \cite{Elmachkouri-Volny-Wu2011}. 
For any $i$ in $\Z^d$ and any $x$ in $\R$, we denote $R_i=\textrm{G}_n(x,i)-\overline{\textrm{G}}_n(x,i)$. Since there exists a measurable function $H$ such that 
$R_i=H(\varepsilon_{i-s};s\in\Z^d)$, we are able to define the physical dependence measure coefficients $(\delta^{(n)}_{i,p})_{i\in\Z^d}$ associated to the random field 
$(R_i)_{i\in\Z^d}$. We recall that $\delta^{(n)}_{i,p}=\|R_i-R_i^{\ast}\|_p$ where $R_i^{\ast}=H(\varepsilon^{\ast}_{i-s};s\in\Z^d)$ and 
$\varepsilon_j^{\ast}=\varepsilon_j\ind{\{j\neq 0\}}+\varepsilon_0^{'}\ind{\{j=0\}}$ for any $j$ in $\Z^d$. In other words, we obtain $R_i^{\ast}$ from $R_i$ 
by just replacing $\varepsilon_0$ by its copy $\varepsilon_0^{'}$. Let $\tau:\Z\to\Z^d$ be a bijection. For any $l\in\Z$, for any $i\in\Z^d$, we denote
\begin{equation}\label{definition_P_i}
P_lR_i=\E(R_i\vert\F_l)-\E(R_i\vert\F_{l-1})
\end{equation}
where $\F_l=\sigma\left(\varepsilon_{\tau(s)};s\leq l\right)$.
\begin{lemma}\label{majoration_P_X} For any $l$ in $\Z$ and any $i$ in $\Z^d$, we have $\|P_lR_i\|_p\leq\delta^{(n)}_{i-\tau(l),p}$.
\end{lemma}
{\em Proof of Lemma $\ref{majoration_P_X}$}. Let $l$ in $\Z$ and $i$ in $\Z^d$ be fixed.
$$
\left\|P_lR_i\right\|_p=\left\|\E(R_i\vert\F_l)-\E(R_i\vert\F_{l-1})\right\|_p=\left\|\E(R_0\vert T^i\F_l)-\E(R_0\vert T^i\F_{l-1})\right\|_p
$$
where $T^i\F_l=\sigma\left(\varepsilon_{\tau(s)-i};s\leq l\right)$.
\begin{align*}
\left\|P_lR_i\right\|_p&=\left\|\E\left(H\left((\varepsilon_{-s})_{s\in\Z^d}\right)\vert T^i\F_l\right)-
\E\left(H\left((\varepsilon_{-s})_{s\in\Z^d\backslash\{i-\tau(l)\}};\varepsilon^{'}_{\tau(l)-i}\right)\vert T^i\F_{l}\right)\right\|_p\\
&\leq \left\|H\left((\varepsilon_{-s})_{s\in\Z^d}\right)-H\left((\varepsilon_{-s})_{s\in\Z^d\backslash\{i-\tau(l)\}};\varepsilon^{'}_{\tau(l)-i}\right)\right\|_p\\
&=\left\|H\left((\varepsilon_{i-\tau(l)-s})_{s\in\Z^d}\right)-H\left((\varepsilon_{i-\tau(l)-s})_{s\in\Z^d\backslash\{i-\tau(l)\}};\varepsilon^{'}_{0}\right)\right\|_p\\
&=\left\|R_{i-\tau(l)}-R_{i-\tau(l)}^{\ast}\right\|_p\\
&=\delta^{(n)}_{i-\tau(l),p}.
\end{align*}
The proof of Lemma $\ref{majoration_P_X}$ is complete.\\ 
\\
For all $i$ in $\Z^d$, $R_i=\sum_{l\in\Z}P_lR_i$. Consequently, $\left\|\sum_{i\in\Lambda_n}a_iR_i\right\|_p=\left\|\sum_{l\in\Z}\sum_{i\in \Lambda_n}a_iP_lR_i\right\|_p$. 
Applying the Burkholder inequality (cf. \cite{Hall--Heyde1980}, page ??) for the martingale difference sequence $\left(\sum_{i\in \Lambda_n}a_iP_lR_i\right)_{l\in\Z}$, we obtain
$$
\left\|\sum_{i\in\Lambda_n}a_iR_i\right\|_p
\leq\left(2p\sum_{l\in\Z}\left\|\sum_{i\in \Lambda_n}a_iP_lR_i\right\|_p^2\right)^{\frac{1}{2}}
\leq\left(2p\sum_{l\in\Z}\left(\sum_{i\in \Lambda_n}\vert a_i\vert\left\|P_lR_i\right\|_p\right)^2\right)^{\frac{1}{2}}
$$
By the Cauchy-Schwarz inequality, we have
$$
\left(\sum_{i\in \Lambda_n}\vert a_i\vert\left\|P_lR_i\right\|_p\right)^2
\leq\sum_{i\in \Lambda_n}a_i^2\left\|P_lR_i\right\|_p\times\sum_{i\in\Lambda_n}\|P_lR_i\|_p
$$
and by Lemma $\ref{majoration_P_X}$, $\sum_{i\in\Z^d}\|P_lR_i\|_p\leq\sum_{j\in\Z^d}\delta^{(n)}_{j,p}$. So, we derive
$$
\left\|\sum_{i\in\Lambda_n}a_iR_i\right\|_p\leq\left(2p\sum_{j\in\Z^d}\delta^{(n)}_{j,p}\sum_{i\in \Lambda_n}a_i^2\sum_{l\in\Z}\left\|P_lR_i\right\|_p\right)^{\frac{1}{2}}.
$$
Applying again Lemma $\ref{majoration_P_X}$, we have $\sum_{l\in\Z}\|P_lR_i\|_p\leq\sum_{j\in\Z^d}\delta^{(n)}_{j,p}$ for any $i$ in $\Z^d$ and finally, we derive
$$
\left\|\sum_{i\in\Lambda_n}a_iR_i\right\|_p\leq\left(2p\sum_{i\in \Lambda_n}a_i^2\right)^{\frac{1}{2}}\sum_{i\in\Z^d}\delta^{(n)}_{i,p}.
$$
Since $\overline{\textrm{G}}_n^{\ast}(x,i)=\E\left(\textrm{G}_n^{\ast}(x,i)\big\vert\F_{n,i}^{\ast}\right)$ where 
$\F_{n,i}^{\ast}=\sigma\left(\varepsilon^{\ast}_{i-s}\,;\,\vert s\vert\leq m_n\right)$ and 
$\left(\textrm{G}_n(x,i)-\overline{\textrm{G}}_n(x,i)\right)^{\ast}=\textrm{G}_n^{\ast}(x,i)-\overline{\textrm{G}}_n^{\ast}(x,i)$, 
we derive $\delta_{i,p}^{(n)}\leq 2\|\textrm{G}_n(x,i)-\textrm{G}_n^{\ast}(x,i)\|_p$. 
Moreover, for any $s$ and $t$ in $\R$, $\vert g_n(s)-g_n(t)\vert\leq\kappa b_n^{-\beta}\vert s-t\vert$. So, we obtain 
\begin{equation}\label{delta_i_p_n_inequality1}
\delta_{i,p}^{(n)}\leq \kappa b_n^{-1-\beta}\delta_{i,p}
\end{equation}
where $\delta_{i,p}=\|X_i-X^{\ast}_i\|_p$.
\begin{lemma}\label{controle_norme_de_la_difference_de_Gnx0_et_Gnx0_barre}
For any $p\geq 2$, any positive integer $n$ and any $x$ in $\R$,
$$
\|\emph{G}_n(x,0)-\overline{\emph{G}}_n(x,0)\|_p\leq\frac{\kappa\sqrt{p}}{b_n^{1+\beta}}\sum_{\vert j\vert>m_n}\delta_{j,p}.
$$
\end{lemma}
{\em Proof of Lemma \ref{controle_norme_de_la_difference_de_Gnx0_et_Gnx0_barre} }. We consider the sequence $(B_n)_{n\geq 0}$ of finite subsets of $\Z^d$ defined by 
$B_0=\{(0,...,0)\}$ and for any $n$ in $\N^{\ast}$, $B_n=\{i\in\Z^d\,;\,\vert i\vert=n\}$. The cardinality of the set $B_n$ is 
$\vert B_n\vert=2d(2n+1)^{d-1}$ for $n\geq 1$. Let $\tau:\N^{\ast}\to\Z^d$ be the bijection defined by $\tau(1)=(0,...,0)$ and 
\begin{itemize}
\item for any $n$ in $\N^{\ast}$, if $l\in\left]a_{n-1},a_n\right]$ then $\tau(l)\in B_n$,
\item for any $n$ in $\N^{\ast}$, if $(p,q)\in\left]a_{n-1},a_n\right]^2$ and $p<q$ then $\tau(p)<_{\textrm{lex}}\tau(q)$
\end{itemize}
where $a_n=\sum_{j=0}^n\vert B_j\vert$ goes to infinity as $n$ goes to infinity. Let $(m_n)_{n\geq 1}$ be the sequence of positive integers defined by (\ref{def_mn}). 
For any $n$ in $\N^{\ast}$, we recall that $\F_{n,0}=\sigma\left(\varepsilon_{-s}\,;\,\vert s\vert\leq m_n\right)$ 
and we consider also the $\sigma$-algebra $\G_{n}:=\sigma\left(\varepsilon_{\tau(p)}\,;\,1\leq p\leq n\right)$. By the definition of the bijection $\tau$, 
for any $n$ in $\N$, $1\leq p\leq a_n$ if and only if $\vert\tau(p)\vert\leq n$. So, we have $\G_{a_{m_n}}=\F_{n,0}$. Consequently, 
$\textrm{G}_n(x,0)-\overline{\textrm{G}}_n(x,0)=\sum_{l> a_{m_n}}D_l$ with 
$D_l=\E\left(\textrm{G}_n(x,0)\vert\G_l\right)-\E\left(\textrm{G}_n(x,0)\vert\G_{l-1}\right)$ for any $l$ in $\Z$. 
Let $p\geq 2$ be fixed. Since $\left(D_l\right)_{l\in\Z}$ is a martingale-difference sequence, applying Burkholder's inequality (cf. \cite{Hall--Heyde1980}, page ??), we derive
$$
\|\textrm{G}_n(x,0)-\overline{\textrm{G}}_n(x,0)\|_p\leq\left(2p\sum_{l>a_{m_n}}\|D_l\|_p^2\right)^{1/2}.
$$
Denoting $\textrm{G}_n^{'}(x,0)=g_n\left(b_n^{-1}\left(x-F\left((\varepsilon_{-s})_{s\in\Z^d\backslash\{-\tau(l)\}};\varepsilon^{'}_{\tau(l)}\right)-\theta_0\right)\right)$, we obtain
$$
\|D_l\|_p=\|\E\left(\textrm{G}_n(x,0)\vert \G_l\right)-\E(\textrm{G}_n^{'}(x,0)\vert\G_l)\|_p\leq \|\textrm{G}_n(x,0)-\textrm{G}_n^{'}(x,0)\|_p.
$$
Keeping in mind that $\vert g_n(s)-g_n(t)\vert\leq \kappa b_n^{-\beta}\vert s-t\vert$ for any $s$ and $t$ in $\R$, we derive
\begin{align*}
\|D_l\|_p&\leq \kappa b_n^{-1-\beta}\left\|F\left((\varepsilon_{-s})_{s\in\Z^d}\right)-F\left((\varepsilon_{-s})_{s\in\Z^d\backslash\{-\tau(l)\}};\varepsilon^{'}_{\tau(l)}\right)\right\|_p\\
&=\kappa b_n^{-1-\beta}\left\|F\left((\varepsilon_{-\tau(l)-s})_{s\in\Z^d}\right)-F\left((\varepsilon_{-\tau(l)-s})_{s\in\Z^d\backslash\{-\tau(l)\}};\varepsilon^{'}_{0}\right)\right\|_p\\
&=\kappa b_n^{-1-\beta}\left\|X_{-\tau(l)}-X_{-\tau(l)}^{\ast}\right\|_p=\kappa b_n^{-1-\beta}\delta_{-\tau(l),p}
\end{align*}
and finally
$$
\|\textrm{G}_n(x,0)-\overline{\textrm{G}}_n(x,0)\|_p\leq\kappa b_n^{-1-\beta}\left(2p\sum_{l>a_{m_n}}\delta_{-\tau(l),p}^2\right)^{1/2}
\leq\kappa b_n^{-1-\beta}\sqrt{p}\sum_{\vert j\vert>m_n}\delta_{j,p}.
$$
The proof of Lemma \ref{controle_norme_de_la_difference_de_Gnx0_et_Gnx0_barre} is complete.\\
\\
Noting that $\delta_{i,p}^{(n)}\leq 2\|\textrm{G}_n(x,0)-\overline{\textrm{G}}_n(x,0)\|_p$ 
and applying Lemma \ref{controle_norme_de_la_difference_de_Gnx0_et_Gnx0_barre}, we derive
\begin{equation}\label{delta_i_p_n_inequality2}
\delta_{i,p}^{(n)}\leq\frac{\kappa\sqrt{p}}{b_n^{1+\beta}}\sum_{\vert j\vert>m_n}\delta_{j,p}.
\end{equation}
Combining ($\ref{delta_i_p_n_inequality1}$) and ($\ref{delta_i_p_n_inequality2}$), we obtain 
$$
\sum_{i\in\Z^d}\delta_{i,p}^{(n)}\leq \kappa\left(\frac{m_n^d\sqrt{p}}{b_n^{1+\beta}}
\sum_{\vert j\vert>m_n}\delta_{j,p}+\frac{1}{b_n^{1+\beta}}\sum_{\vert j\vert>m_n}\delta_{j,p}\right).
$$
The proof of Lemma \ref{moment_inequality} is complete.\\
\\
{\em Proof of Lemma $\ref{lemme-technique-barre}$}. Let $s$ and $t$ be fixed in $\R$. We have 
\begin{align*}
\big\vert\E\left(\overline{\textrm{G}}_n(s,0)\overline{\textrm{G}}_n(t,0)\right)-\E\left(\textrm{G}_n(s,0)\textrm{G}_n(t,0)\right)\big\vert
&\leq \|\textrm{G}_n(s,0)\|_2\|\textrm{G}_n(t,0)-\overline{\textrm{G}}_n(t,0)\|_2\\
&\qquad\qquad+\|\overline{\textrm{G}}_n(t,0)\|_2\|\textrm{G}_n(s,0)-\overline{\textrm{G}}_n(s,0)\|_2.
\end{align*}
Using (\ref{variance_gn}) and applying Lemma \ref{controle_norme_de_la_difference_de_Gnx0_et_Gnx0_barre}, we derive
$$
\big\vert\E\left(\overline{\textrm{G}}_n(s,0)\overline{\textrm{G}}_n(t,0)\right)-\E\left(\textrm{G}_n(s,0)\textrm{G}_n(t,0)\right)\big\vert
\leq\frac{\kappa}{b_n^{2\beta+\frac{1}{2}}}\sum_{\vert j\vert>m_n}\delta_j.
$$
Since $b_n^2\vert\E(Z_0(s)Z_0(t))-\E(\overline{Z}_0(s)\overline{Z}_0(t)\vert=\vert\E\left(\textrm{G}_n(s,0)\textrm{G}_n(t,0)\right)-
\E\left(\overline{\textrm{G}}_n(s,0)\overline{\textrm{G}}_n(t,0)\right)\vert$, we obtain
\begin{equation}\label{ecart_Z_0_Z_0_barre}
m_n^db_n^{2\beta+1}\vert\E(Z_0(s)Z_0(t))-\E(\overline{Z}_0(s)\overline{Z}_0(t)\vert
\leq\frac{\kappa}{(m_n^db_n)^{3/2}}\sum_{\vert j\vert>m_n}\vert j\vert^{\frac{5d}{2}}\,\delta_j.
\end{equation}
Combining (\ref{convergence_esperance_Zs_carre_bis}), (\ref{ecart_Z_0_Z_0_barre}) and Lemma \ref{mn}, for any $z$ in $\R$, we obtain 
\begin{equation}\label{limite_variance_Z0_barre}
b_n^{2\beta+1}\E(\overline{Z}_0^2(z))\converge{n}{+\infty}{ }\sigma^2(z).
\end{equation}
Similarly, combining (\ref{limite_esperance_Zx_Zy}), (\ref{ecart_Z_0_Z_0_barre}) and Lemma \ref{mn}, we derive
\begin{equation}\label{limite_esperance_Z0xZ0y_barre}
\vert\E(\overline{Z}_0(x)\overline{Z}_0(y))\vert=o(m_n^{-d}b_n^{-2\beta-1}).
\end{equation}
In the other part, we have
\begin{equation}\label{decomposition_esperance_delta_0_carre_barre}
\E(\overline{\Delta}_0^2)=\lambda_1^2\E(\overline{Z}_0^2(x))+\lambda_2^2\E(\overline{Z}_0^2(y))+2\lambda_1\lambda_2\E(\overline{Z}_0(x)\overline{Z}_0(y)).
\end{equation}
Combining (\ref{limite_variance_Z0_barre}), (\ref{limite_esperance_Z0xZ0y_barre}) and (\ref{decomposition_esperance_delta_0_carre_barre}), we obtain the convergence 
of $b_n^{2\beta+1}\E(\overline{\Delta}_0^2)$ to $\eta$ as $n$ goes to infinity.\\
Let $i\neq 0$ be fixed in $\Z^d$ and let $s$ and $t$ be fixed in $\R$. We have
\begin{equation}\label{borne_Z_0_Z_i}
m_n^db_n^{2\beta+1}\E\vert\overline{Z}_0(s)\overline{Z}_i(t)\vert\leq m_n^db_n^{2\beta-1}\left(\E\big\vert \overline{\textrm{G}}_n(s,0)\overline{\textrm{G}}_n(t,i)\big\vert
+3\E\big\vert\textrm{G}_n(s,0)\big\vert\,\E\big\vert \textrm{G}_n(t,0)\big\vert\right).
\end{equation}
Since $\vert\vert \alpha\vert-\vert \beta\vert\vert\leq\vert \alpha-\beta\vert$ for any $(\alpha,\beta)$ in $\R^2$ and applying the Cauchy-Schwarz inequality, we obtain 
\begin{align*}
\big\vert\E\vert\overline{\textrm{G}}_n(s,0)\overline{\textrm{G}}_n(t,i)\vert-\E\vert \textrm{G}_n(s,0)\textrm{G}_n(t,i)\vert\big\vert
&\leq\|\overline{\textrm{G}}_n(s,0)\|_2\|\overline{\textrm{G}}_n(t,0)-\textrm{G}_n(t,0)\|_2\\
&\qquad\qquad+\|\textrm{G}_n(t,0)\|_2\|\overline{\textrm{G}}_n(s,0)-\textrm{G}_n(s,0)\|_2
\end{align*}
Using again (\ref{variance_gn}) and applying Lemma \ref{controle_norme_de_la_difference_de_Gnx0_et_Gnx0_barre}, we derive
\begin{equation}\label{ecart_overline_GnO_Gni_et_Gn0_Gni_barre}
m_n^db_n^{2\beta-1}\big\vert\E\vert\overline{\textrm{G}}_n(s,0)\overline{\textrm{G}}_n(t,i)\vert-\E\vert\textrm{G}_n(s,0)\textrm{G}_n(t,i)\vert\big\vert
\leq\frac{\kappa}{(m_n^db_n)^{3/2}}\sum_{\vert j\vert>m_n}\vert j\vert^{\frac{5d}{2}}\delta_j.
\end{equation}
Combining (\ref{covariance_de_gn}), (\ref{ecart_overline_GnO_Gni_et_Gn0_Gni_barre}) and Lemma \ref{mn}, we obtain 
\begin{equation}\label{limite_mndbn_covariance_Gn0}
\sup_{i\in\Z^d\backslash\{0\}}\E\vert\overline{\textrm{G}}_n(s,0)\overline{\textrm{G}}_n(t,i)\vert=o(m_n^{-d}b_n^{-2\beta+1}).
\end{equation}
Moreover, using ($\ref{esperance_gn}$) and applying again Lemma \ref{mn}, we have 
\begin{equation}\label{limite_mndbn_covariance_Gn0_bis}
\E\big\vert\textrm{G}_n(s,0)\big\vert\,\E\big\vert \textrm{G}_n(t,0)\big\vert=o(m_n^{-d}b_n^{-2\beta+1}).
\end{equation}
Combining (\ref{borne_Z_0_Z_i}), (\ref{limite_mndbn_covariance_Gn0}) and (\ref{limite_mndbn_covariance_Gn0_bis}), we derive
$$
\sup_{i\in\Z^d\backslash\{0\}}\E\vert\overline{Z}_0(s)\overline{Z}_i(t)\vert=o(m_n^{-d}b_n^{-2\beta-1}).
$$
Since 
$\E\vert\overline{\Delta}_0\overline{\Delta}_i\vert\leq\lambda_1^2\E\vert\overline{Z}_0(x)\overline{Z}_i(x)\vert
+2\lambda_1\lambda_2\E\vert \overline{Z}_0(x)\overline{Z}_i(y)\vert+\lambda_2^2\E\vert \overline{Z}_0(y)\overline{Z}_i(y)\vert$, we have also 
$\sup_{i\in\Z^d\backslash\{0\}}\E\vert\overline{\Delta}_0\overline{\Delta}_i\vert=o(m_n^{-d}b_n^{-2\beta-1})$. The proof of Lemma $\ref{lemme-technique-barre}$ is complete.
% ------------------------------------------------------------------------
\bibliographystyle{plain}
\bibliography{xbib}

\begin{thebibliography}{10}

\bibitem{Bosq-Merlevede-Peligrad1999}
D.~Bosq, Merlev\`ede F., and M.~Peligrad.
\newblock Asymptotic normality for density kernel estimators in discrete and
  continuous time.
\newblock {\em J. Multivariate Anal.}, 68:78--95, 1999.

\bibitem{Casti1985}
John~L. Casti.
\newblock {\em Nonlinear system theory}, volume 175 of {\em Mathematics in
  Science and Engineering}.
\newblock Academic Press Inc., Orlando, FL, 1985.

\bibitem{Dedecker1998}
J.~Dedecker.
\newblock A central limit theorem for stationary random fields.
\newblock {\em Probab. Theory Relat. Fields}, 110:397--426, 1998.

\bibitem{Elmachkouri2011}
M.~El~Machkouri.
\newblock Asymptotic normality for the parzen-rosenblatt density estimator for
  strongly mixing random fields.
\newblock {\em Statistical Inference for Stochastic Processes}, 14(1):73--84,
  2011.

\bibitem{Elmachkouri--Stoica2010}
M.~El~Machkouri and R.~Stoica.
\newblock Asymptotic normality of kernel estimates in a regression model for
  random fields.
\newblock {\em J. Nonparametr. Stat.}, 22(8):955--971, 2010.

\bibitem{Elmachkouri-Volny-Wu2011}
M.~El~Machkouri, D.~Voln\`y, and W.~B. Wu.
\newblock A central limit theorem for stationary random fields.
\newblock arXiv:1109.0838v1, 2011.

\bibitem{Fan1991b}
J.~Fan.
\newblock Asymptotic normality for deconvolution kernel density estimators.
\newblock {\em Sankhya Ser. A}, 53(1):97--110, 1991.

\bibitem{Fan1991a}
J.~Fan.
\newblock On the optimal rates of convergence for nonparametric deconvolution
  problems.
\newblock {\em Ann. Statist.}, 19(3):1257--1272, 1991.

\bibitem{Hall--Heyde1980}
P.~Hall and C.~C. Heyde.
\newblock {\em Martingale limit theory and its application}.
\newblock Academic {P}ress, New {Y}ork, 1980.

\bibitem{Li2008}
Z.~Li.
\newblock Asymptotic normality for deconvolution kernel density estimators from
  random fields.
\newblock arXiv:0810.3121v1, 2008.

\bibitem{Lindeberg1922}
J.~W. Lindeberg.
\newblock Eine neue {H}erleitung des {E}xponentialgezetzes in der
  {W}ahrscheinlichkeitsrechnung.
\newblock {\em Mathematische {Z}eitschrift}, 15:211--225, 1922.

\bibitem{Masry1993b}
E.~Masry.
\newblock Asymptotic normality for deconvolution estimators of multivariate
  densities of stationary processes.
\newblock {\em J. Multivariate Anal.}, 44(1):47--68, 1993.

\bibitem{Masry1993a}
E.~Masry.
\newblock Strong consistency and rates for deconvolution of multivariate
  densities of stationary processes.
\newblock {\em Stochastic Process. Appl.}, 47(1):53--74, 1993.

\bibitem{McLeish}
D.~L. McLeish.
\newblock A maximal inequality and dependent strong laws.
\newblock {\em Ann. Probab.}, 3(5):829--839, 1975.

\bibitem{Rio1993}
E.~Rio.
\newblock Covariance inequalities for strongly mixing processes.
\newblock {\em Ann. Inst. H. Poincar\'e Probab. Statist.}, 29(4):587--597,
  1993.

\bibitem{Ros}
M.~Rosenblatt.
\newblock A central limit theorem and a strong mixing condition.
\newblock {\em Proc. {N}at. {A}cad. {S}ci. USA}, 42:43--47, 1956.

\bibitem{Rugh1981}
W.~J. Rugh.
\newblock {\em Nonlinear system theory}.
\newblock Johns Hopkins Series in Information Sciences and Systems. Johns
  Hopkins University Press, Baltimore, Md., 1981.

\bibitem{Stefanski--Carroll1990}
L.~Stefanski and R.~J. Carroll.
\newblock Deconvoluting kernel density estimators.
\newblock {\em Statistics}, 21(2):169--184, 1990.

\bibitem{Tran1990}
L.T. Tran.
\newblock Kernel density estimation on random fields.
\newblock {\em J. Multivariate Anal.}, 34:37--53, 1990.

\bibitem{Wu2005}
W.~B. Wu.
\newblock Nonlinear system theory: another look at dependence.
\newblock {\em Proc. Natl. Acad. Sci. USA}, 102(40):14150--14154 (electronic),
  2005.

\bibitem{Wu2011}
W.~B. Wu.
\newblock Asymptotic theory for stationary processes.
\newblock {\em Statistics and Its Interface}, 0:1--20, 2011.

\end{thebibliography}
% ------------------------------------------------------------------------
\end{document}